\documentclass[12pt]{article}

\usepackage{amsmath,amssymb,amsthm}
\usepackage{verbatim}
\usepackage{color}
\usepackage{graphics}

\newtheorem{theorem}{Theorem}[section]
\newtheorem{thm}[theorem]{Theorem}
\newtheorem{corollary}[theorem]{Corollary}

\newtheorem{lemma}[theorem]{Lemma}
\newtheorem{proposition}[theorem]{Proposition}

\newtheorem{definition}[theorem]{Definition}
\newtheorem{remark}[theorem]{Remark}
\newtheorem{example}[theorem]{Example}
\numberwithin{equation}{section}

\def\be{\begin{equation}}
\def\ee{\end{equation}}
\def\bes{\begin{equation*}}
\def\ees{\end{equation*}}

\def\vp{{\varphi}}
\def\q{\quad}
\def\qq{{\qquad}}

\def\pd{\partial}

\def\wt{\widetilde}

\def\eps{\varepsilon} 
\def\al{\alpha}
\def\phi{\varphi}
\def\lam{{\lambda}}

\def\ol{\overline}
\def\fract{\textstyle \frac}
\def\Gam{\Gamma} \def\gam{\gamma}
\def\th{\theta}

\def\proof{{\medskip\noindent {\bf Proof. }}}

\def\qed{{\hfill $\square$ \bigskip}}

\def\supp{{\mathop {{\rm supp\, }}}}

\def\sD {{\cal D}} \def\sE {{\cal E}} \def\sF {{\cal F}}
  
  \def\sL {{\cal L}}
\def\sM {{\cal M}} \def\sN {{\cal N}}

  \def\sX {{\cal X}}

 \def\bE {{\mathbb E}}

 \def\bN {{\mathbb N}} 
\def\bP {{\mathbb P}}  \def\bR {{\mathbb R}}

\def\nn{\nonumber}

\def\half{{\textstyle \frac12}}
\def\ignore#1{}  

\def\ms{\medskip}
\def\sms{\smallskip}

\def\sm{\smallskip\noindent}

\def\grad{\nabla}

\definecolor{dred}{rgb}{0.8, 0.0, 0.0}

\def\indicator{{\mathchoice {1\mskip-4mu\mathrm l}%
{1\mskip-4mu\mathrm l}{1\mskip-4.5mu\mathrm l}%
{1\mskip-5mu\mathrm l}}}

\def\UHK{\mathrm{UHK}}
\def\HK{\mathrm{HK}}
\def\CSA{\mathrm{CSA}}
\def\CSD{\mathrm{CSD}}
\def\CS{\mathrm{CS}}
\def\PI{\mathrm{PI}}
\def\FK{\mathrm{FK}}
\def\tFK{\mathit{FK}}
\def\FK{\mathrm{FK}}
\def\VD{\mathrm{VD}}
\def\Psiinv{{\Psi^{-1}}}

\def\esssup{{\mathop{\rm ess \; sup \, }}}

\def\rho{{\varrho}}

\begin{document}

\title{\bf Energy inequalities for cutoff functions and some  applications}

\author{Sebastian Andres\footnote{Research partially supported 
by NSERC (Canada)}, 
Martin T. Barlow\footnote{Research partially supported by 
NSERC (Canada) and by Trinity
College, Cambridge}}

\date{}
\maketitle

\begin{abstract}
Let $(\sX,d,m)$ be a metric measure space with a local regular Dirichlet form. 
We establish necessary and sufficient conditions for upper heat kernel  
bounds with  sub-diffusive space-time exponent  to hold. This characterization is 
stable under rough isometries, that is it is preserved under bounded perturbations 
of the Dirichlet form.
Further, we give a criterion for stochastic completeness in terms of a Sobolev 
inequality for cutoff functions. As an example we show that this criterion applies 
to an anomalous diffusion on a geodesically incomplete fractal space, where the 
well-established criterion in terms of volume growth fails.

\end{abstract}

\section{Introduction} \label{sec:intro}

Let $(\sX,d)$ be a locally compact metric space and let $m$ be
a positive Radon measure on $\sX$ with $\supp[m]=\sX$. 
We will refer to such a triple $(\sX, d, m)$ as a \emph{metric measure space},
and denote by 
$\langle .,.\rangle$ the inner product in $L^2(\sX,m)$.
We consider a regular, strongly local Dirichlet form $(\sE, \sF)$
on $L^2(\sX,m)$ (see \cite{FOT}). Let 
$\sL$ be the (negative definite) generator of $\sE$; this is a 
self-adjoint operator in $L^2(\sX,m)$ such that
$$ \sE (f,g) = -\langle \sL f,g \rangle \q \hbox{ for all 
$f \in \mathcal{D}(\sL)$, $g \in \sF$}, $$
and let $\{P_t\}_{t\geq 0}$  be the associated semigroup.
If $P_t$ has a density $p_t(x,y)$ with respect to $m$
then after some regularization we 
call this the {\em heat kernel} on the 
\emph{metric measure Dirichlet space} (or \emph{MMD space})
$(\sX, d, m, \sE)$. 
Our main interest is in upper bounds on $p_t(x,y)$.
Write $B(x,r)$ for balls in $(\sX,d)$ and set
\be
  V(x,r) = m(B(x,r)). 
\ee 
Most familiar are Gaussian upper bounds of the form
\be \label{e:gub}
 p_t(x,y) \le \frac{c_1}{V(x,t^{1/2})} \exp\Big( -c_2 \frac{d(x,y)^2}{t}\Big); 
\ee
these arise (with lower bounds of the same form but with different constants) 
in the case of uniformly elliptic divergence form PDE, 
and manifolds with Ricci curvature bounded uniformly below  -- see \cite{Ar,LY}. 

If \eqref{e:gub} holds we will say $(\sX,\sE)$ satisfies the condition
$\UHK(2)$; if in addition Gaussian lower bounds hold we say
$\HK(2)$ holds. One can ask for characterizations of these bounds, and
in particular for characterizations which are {\em stable}, 
that is that are  preserved under bounded perturbation of the
Dirichlet form. More precisely, a property (P) of $(\sX, \sE)$ is
stable if when $(\sE_i, \sF)$ are
two Dirichlet forms on $L^2(\sX,m)$ with
$$  C^{-1} \sE_1(f,f) \le \sE_2(f,f) \le C \sE_1(f,f), \q f \in \sF, $$
then $(P)$ holds for $(\sX,\sE_1)$ if and only if it holds for $(\sX,\sE_2)$.
In the manifold case 
stability for $\HK(2)$ was proved in \cite{Gr0, SC1} by showing that 
these Gaussian bounds are equivalent to {\em volume doubling} (denoted VD)
plus a family of Poincar\'e inequalities -- see below for the precise
definitions. If VD holds then stability for $\UHK(2)$ is a consequence 
of the results of \cite{Gr3}, where it is shown that $\UHK(2)$ 
is equivalent to a Faber Krahn inequality $\FK(2)$, which controls
the smallest eigenvalue of domains in $\sX$. 

\sms
The Gaussian bounds \eqref{e:gub} arise due to the standard space-time scaling
relation $t = r^2$. More general possibilities can arise; for various exact
fractals (see \cite{Ba1}) one can have $V(x,r) \asymp r^\al$ and a space
time scaling of $t=r^\beta$, where 
$\al \in [1,\infty)$ and 
 $\beta \in [2,1+ \al]$;
the case when $\beta \neq 2 $ is called {\em anomalous diffusion.}
 Since we wish to
be able to consider spaces with different local and global structure, 
we introduce a more general space-time scaling function $\Psi$. 
Let $\beta_L \ge 2$, $\beta \ge 2$, and set 
\be \label{e:Psidef}
\Psi(r)= \Psi_{\beta_L, \beta}(r) =
\begin{cases}  
  r^{\beta_L} & \hbox{ if  } 0\le r \le 1, \\
   r^\beta   & \hbox{ if  } r > 1.
\end{cases}  
\ee
We will write $\UHK(\Psi)$ for the heat kernel upper bounds associated
with $\Psi$ -- see Definition \ref{D:uhk} below for their precise form.
Our main theorem is a stable characterization of  $\UHK(\Psi)$, in terms
of Faber-Krahn inequality $\FK(\Psi)$, and a new condition 
denoted $\CSA(\Psi)$, which controls
the energy of cutoff functions in annuli.

\ms
To state our results precisely, we need a number of further
definitions.

Since $\sE$ is regular, each function $f\in \sF$ admits a
quasi-continuous version $\tilde f$ (see Theorem 2.1.3 in
\cite{FOT}). Throughout the paper, we will abuse notation and take
the quasi-continuous version of $f$ without writing $\tilde f$. 
Another consequence of regularity is that $\sE(f,g)$ can be
written in terms of a signed measure $\Gamma(f,g)$ as
$$ \sE(f,g)=\int_\sX d\Gamma(f,g). $$
For any essentially bounded $f\in \sF$, $\Gamma(f,f)$ is the unique 
Borel measure on $\sX$ (called the energy measure) on $\sX$ satisfying
$$
\int_\sX g \, d\Gamma(f,f)=2 \sE(f, fg)-\sE(f^2,g)
$$
for all essentially bounded $g\in \sF$; 
$\Gamma(f,g)$ is then defined by polarization.

\sm {\bf Example.} (Davies \cite{D}). Let $(M,d)$ be a manifold with 
Riemannian volume measure $\mu$, and $\sE(f,f)=\int |\grad f|^2 d\mu$. 
Let $\sigma>0$, and $dm=\sigma^2 \, d\mu$. Then
$$  \sL f = \sigma^{-2}  \grad ( \sigma^2 \grad f ), $$
and $d\Gam(f,f) = |\grad f|^2 \sigma^2 d\mu$.

\sms
For later use we collect from \cite[Section 3.2]{FOT} 
some properties of the energy measure.
\begin{enumerate}
\item[i)] {\bf Locality.} 
For all functions $f,g\in \sF$ and all measurable 
sets $G \subset \sX$ on which $f$ is constant
$$ \indicator_G \, d\Gamma(f,g)=0. $$
\item[ii)] {\bf Leibniz and chain rules.} 
For $f,g \in \sF$ essentially bounded and $\phi \in C^1(\bR)$,
\begin{align*}
d\Gamma(f g,h)&=f \, d\Gamma(g,h)+g\, d\Gamma(f,h),  \\
d\Gamma(\phi(f),g)&=\phi'(f) \, d\Gamma(f,g).
\end{align*} 
\end{enumerate}

We note also the following result of
Le Jan \cite[Proposition 1.5.5(b)]{LJ} -- see also \cite{Mos}, p.\ 389
for a simple proof.

\begin{lemma} \label{L:sgam}
Let $\sX$ be a MMD space. 
Suppose that $(\sE_{i}, \sF), i=1,2$, are 
strongly local regular Dirichlet forms  that satisfy 
\be
 C^{-1} \sE_1(f,f) \le \sE_2(f,f) \le C \sE_1(f,f), \q f \in \sF.
\ee
Then their energy measures $\Gam^{(i)}$ satisfy
\be C^{-1}  d\Gam^{(1)}(f,f) \le d\Gam^{(2)}(f,f)\le C d\Gam^{(1)}(f,f),
\hbox{ for all } f\in \sF.
\ee
\end{lemma}

\ms

We now introduce a number of conditions which the space
$\sX$ and Dirichlet form $\sE$ may or may not satisfy.

\begin{definition} {\rm
We say that $(\sX,d,m)$ satisfies {\sl volume doubling} (VD) if
there exists a constant $C_D$ such that
for every $x \in \sX$, $r > 0$,
\be \label{e:vd}
 V(x,2r) \le C_D V(x,r). 
\ee
} \end{definition}

We next introduce the Faber-Krahn inequality: see
\cite{GT}, Section 3.3 for more details. 
For any open set $D \subset \sX$, $\sF_D$ is defined to be the 
closure in $\sF$ of the set of all functions in $\sF$ that 
are compactly supported in $D$. 
For $D \subset \sX$
we write $\lam_1(D)$ for the smallest (Dirichlet) eigenvalue of
$\sL$ on $D$; this can be defined by the variational formula
\be \label{e:lam1}
 \lam_1(D) 
= \inf \Big\{ \frac{ \sE(f,f)}{||f||_2^2}: f \in \sF_D, f \neq 0 \Big\}.
\ee

\begin{definition}
{\rm 
The MMD space $(\sX, \sE)$ satisfies the {\em Faber-Krahn inequality}
$\FK(\Psi)$ if there exists a constant $C_F$ and $\nu>0$ such that 
for any ball $B=B(x,r)$ and open set $D \subset B$,
\be \label{e:fki}
 \lam_1 (D) \ge \frac{C_F}{\Psi(r)} (m(B)/m(D))^{\nu}. 
\ee
} \end{definition}

We remark that the value of $\nu$ turns out to be unimportant.

\begin{definition} {\rm 
We say that the {\em Poincar\'e inequality} 
$\PI(\Psi)$ holds if there exists a constant $C_P$ such that
for all balls $B=B(x,r)$  and $f \in \sF$, 
$$ \inf_{a \in \bR} \int_B (f-a)^2 dm
= \int_B (f-\ol f_B)^2 dm \le C_P \Psi(r)\int_B d\Gam(f,f). $$
Here $\ol f_B$ is the mean of $f$ on $B$.
} \end{definition}

Associated with the Dirichlet form $(\sE, \sF)$ and semigroup
$(P_t)$ is a Hunt
process $X=(X_t, t \ge 0, \bP^x, x \in \sX-\sN)$. 
Here $\sN$ is `properly exceptional': $m(\sN)=0$ and 
$\bP^x( X_t \in \sN \hbox{ for some } t>0 )=0$ for all $x \in \sX-\sN$ --
see \cite[p.\ 134]{FOT}. This Hunt process is unique up to a 
properly exceptional set -- see \cite[Theorem 4.2.7]{FOT}.
We fix $X$ and $\sN$, and write
\be
  \sX_0 = \sX - \sN.
\ee
While the semigroup $(P_t)$ associated with $\sE$ is defined on $L^2$,
a more precise version, with better regularity properties, can be obtained
if we set, for bounded Borel $f$,
$$ P_t f(x) = \bE^x f(X_t), \q x \in \sX_0. $$
The heat kernel associated with $(P_t)$ (if it exists) is a 
measurable function 
$p_t(x,y): (0,\infty) \times \sX_0 \times \sX_0 \to (0,\infty)$ such that
\begin{align} \label{e:hkdef}
\bE^x f(X_t) &= P_tf(x) = \int p_t(x,y) f(y) \, m(dy), \, x \in \sX_0, 
f \in L^\infty(\sX), \\
p_t(x,y) &= p_t(y,x), \, \hbox{ for all } t>0,\, x ,y \in \sX_0, \\
\label{e:ck}
p_{s+t}(x,z) &= \int p_s(x,y) p_t(y,z) \, m(dy), \hbox{ for all } s>0, t>0,
\,\, x,z \in \sX_0. 
\end{align}
While \eqref{e:hkdef} only defines $p_t(x, \cdot)$ $m$-a.e.,
using the Chapman-Kolmogorov equation \eqref{e:ck} one can 
regularise $p_t(x,y)$ so that \eqref{e:hkdef}--\eqref{e:ck} hold on all
of $\sX_0$. For more details see \cite{GT}. 

Define the function 
\be
 \Phi(R,t) = \sup_{s>0} \Big( \frac{R}{s} - \frac{t}{\Psi(s)} \Big).
\ee
The following lemma summarises some properties of this function -- see
Section 3.3 of \cite{GT} and in particular Example 3.18.

\begin{lemma}\label{L:Phi}
$\Phi(R,t)$ is non-negative, increasing in $R$ and decreasing in $t$. We have
\be 
 \Phi(R,t) \asymp
\begin{cases}
 \Big( \frac{R^{\beta_L}}{t} \Big)^{1/(\beta_L-1)}, & \hbox{ if } t  \le R, \\
  \Big( \frac{R^\beta}{t} \Big)^{1/(\beta-1)}, & \hbox{ if } t \ge R.
\end{cases}
\ee
Further $\Phi(R,\Psi(R)) \le \beta_2^{-1/(\beta_2-1)}$, where $\beta_2= \beta_L \vee \beta$.
\end{lemma}
We define $\Psiinv$ to be the inverse of $\Psi$, so that
$$ \Psiinv(s) = s^{1/\beta_L} 1_{( s\le 1)} + s^{1/\beta} 1_{(s>1)}. $$

\begin{definition} \label{D:uhk}
{\rm
We say $(p_t)$ satisfies $\UHK(\Psi)$ if there exists a properly exceptional
set $\sN_0$ and constants $c_1$, $c_2$ such that 
\be \label{e:uhkdef}
 p_t(x,y) \le V(x, \Psi^{-1}( c_1 t) )^{-1}  \exp( - \Phi(c_2 d(x,y),t) )
\ee
for all $t>0$ and for all $x,y \in \sX-\sN_0$. 
If a similar lower bound (with different constants
$c_i$) also holds then we say that $\HK(\Psi)$ holds.
} \end{definition}

When $\Psi(r)=r^\beta$ we will write $\PI(\beta)$ etc. for 
the condition $\PI(\Psi)$.

\ms
As explained above, we wish to find a stable characterization of the
heat kernel bounds $\UHK(\Psi)$. 
In view of Lemma \ref{L:sgam}, the characterizations of $\HK(2)$ and
$\UHK(2)$ in terms of Faber-Krahn and Poincar\'e inequalities are stable.
It is easy to see that the natural
generalization of these to more general $\Psi$ fails.
Let $\Psi_2 \ge \Psi_1$ with $\Psi_2(r)/\Psi_1(r) \to \infty$,
and let $\sX$ be an unbounded space satisfying $\HK(\Psi_1)$. Then
$\sX$ also satisfies $\FK(\Psi_1)$ and  $\PI(\Psi_1)$, and so by
the monotonicity of these conditions in $\Psi$, it is immediate that 
$\sX$ satisfies $\FK(\Psi_2)$ and  $\PI(\Psi_2)$. However, 
it is straightforward to check that 
$\HK(\Psi_1)$ and $\UHK(\Psi_2)$ cannot both hold.
At a more fundamental level, the conditions $\PI(\Psi)$ and $\FK(\Psi)$
ensure that the heat equation homogenises over a ball of radius $R$
in time at most $\Psi(R)$, but do not exclude the possibility that
this might occur more quickly. To `capture' $\HK(\Psi)$ one needs 
a condition which gives an upper bound on the rate at which heat,
or the diffusion $X$, can move on the space $\sX$. 
Such a condition was found in \cite{BB3, BBK}, which 
gave a stable characterization of $\HK(\Psi)$. 

\begin{definition}
{\rm 
Let $U \subset V$ be open sets in $\sX$ with
$U \subset \ol U \subset V$.
We say a continuous function $\vp$ is a {\sl cutoff function for $U \subset V$}
if $\vp  = 1$ on $U$ and $\vp=0$ on $V^c$.
} \end{definition}

\begin{definition} [\bf Condition ($\CS(\Psi$)]
{\rm (See \cite{BB3, BBK}.) 
We say that condition $\CS(\Psi)$ holds if 
there exist constants $c_1$ and $\th\in (0,1]$ such that the following
holds. For every ball $B(x_0,r)$
there exists a cutoff function $\vp$ with $\vp=1$ on $B(x_0, r/2)$
and $\vp=0$ on $B(x_0, r)^c$, with the following properties. \\
(1) $\vp$ is H\"older continuous of order $\th$. \\
(2) If $0< s \le r$ and $f \in \sF$ then
\be \label{e:cs3}
 \int_{B(y,s)} f^2 \, d\Gamma(\vp,\vp) 
\le c_1 (s/r)^{2 \theta} \Big( \int_{B(y,2s)} d\Gamma(f,f) 
+\Psi(s)^{-1} \int_{B(y,2s)} f^2 dm  \Big).
\ee 
} \end{definition}

`CS' here refers to `cutoff Sobolev'; this condition ensures
the existence of a large class of cutoff functions with low energy.
The main theorem of \cite{BB3, BBK} is that
$\HK(\Psi)$ is equivalent to VD $+ \PI(\Psi)+ \CS(\Psi)$.
While the condition $\CS(\Psi)$ is hard to verify, it is stable. Further,
this stability allows estimates on (for example) the heat kernel on the
Sierpinski carpet to be transferred to manifolds, graphs, or domains in $\bR^d$
which are roughly isometric to the Sierpinski carpet. For rough
isometries see \cite{Kan}, and a for more detailed discussion of this point
see \cite[Section 5]{BBK}.

We now introduce a simplication of the condition $\CS(\Psi)$, which controls
the energy of cutoff functions in annuli.

\begin{definition} {\rm 
Let $D_0$, $D_1$ be open subsets of $\sX$ with 
$D_0 \subset \ol D_0 \subset D_1$, and let $U= D_1 - \ol D_0$.
We say that condition $\CSD(D_0, D_1, \th)$ holds if 
there exists a cutoff function $\vp$ for $D_0 \subset D_1$ 
such that if $f \in \sF$ then, 
\be \label{e:csa1}
 \int_U f^2 \, d\Gamma(\vp,\vp) 
\le {\frac18} \int_U \vp^2 d\Gamma(f,f) + \th  \int_U  f^2  dm. 
\ee 
} \end{definition}

\begin{definition}[\bf Condition ($\CSA(\Psi)$]
{\rm We say that condition $\CSA(\Psi)$ holds if there exists a 
constant $C_S$ such that for every 
$x \in \sX$, $R>0$, $r>0$ 
the condition $\CSD( B(x, R), B(x, R+r), C_S \Psi(r)^{-1} )$ holds.
} \end{definition}

\begin{remark} 
{\rm 
1. If VD holds then $\CS(\Psi)$ implies $\CSA(\Psi)$ 
-- see Lemma \ref{L:cs2csa}.  \\
2. Note that $\CSA(\Psi)$ does not require the 
H\"older continuity of the cutoff function. \\
3. It is essential for the use of \eqref{e:csa1} in Lemma  \ref{L:ci} 
that the constant in front of the first term on the right hand side is less 
than $\fract14$. However, as we will see in Section \ref{s:csa}, the
inequality $\CSA(\Psi)$ has a `self-improving' property, 
which enables one to 
alter the weights of the two terms on the right-hand side. \\
4. It is easy to see by using linear cutoff functions that $\CSA(2)$ 
holds on any manifold.\\
5. The bound \eqref{e:csa1} is not symmetric between $\vp$ and
$1-\vp$, but very often we will just use the fact that $\vp \le 1$
in the first term on the right hand side. \\
6.  In view of Lemma \ref{L:sgam} and the results of Section 5, 
the condition $\CSA(\Psi)$ is stable -- see Corollary \ref{C:csastab}. \\
7. See \cite{Bas} for the use of an inequality similar to 
$\CS(\Psi)$ to prove stability of the elliptic Harnack inequality 
for a class of graphs.  
}
\end{remark}

\sms
Our first main theorem is the following characterization
of $\UHK(\Psi)$. 

\begin{thm} \label{T:main} 
Assume that $\sX$ satisfies VD and is unbounded in the metric $d$.
The following are equivalent: \\
(1) $\FK(\Psi)$ and $\CSA(\Psi)$. \\
(2) $\UHK(\Psi)$. \\
\end{thm}

\begin{remark}\label{R:main}
{\rm 
(1) See \cite{GH} for several other conditions equivalent
to $\UHK(\beta)$. 
Note however that unlike (1) above, none of these were known to
be stable under bounded perturbation of the Dirichlet form $\sE$. \\
(2) The reader may wonder if while $\CSA(\Psi)$ is sufficient for upper bounds,
one needs the stronger $\CS(\Psi)$ to obtain lower bounds as well.
However, we expect that $\HK(\Psi)$ is equivalent
to $\VD+ \PI(\Psi) + \CSA(\Psi)$. In fact, $\VD + \PI(\Psi)$ is enough
to give $\FK(\Psi)$, so using Theorem \ref{T:main} 
one obtains $\UHK(\Psi)$.  
Given this, the methods of Stroock and Saloff-Coste \cite{SCS},
and Fabes-Stroock \cite{FS} should then lead to a matching lower bound. 
} \end{remark}

The main theorem of \cite{BB3, BBK} was proved using Moser's method 
\cite{Mo1}. To prove the implication $(1) \Rightarrow (2)$ in Theorem \ref{T:main}
we first show in Proposition \ref{P:DG} that $\CSA(\Psi)$ gives a generalization 
of the `Davies-Gaffney' bound of \cite{D}.
Next, we use techniques developed in
\cite{Gr0, CG} to prove a mean value inequality for caloric functions
(i.e. solutions of the heat equation), which leads to the pointwise bounds
$\UHK(\Psi)$. For the easier implication $(2) \Rightarrow (1)$
we use the method of \cite{BBK}, but since $\CSA(\Psi)$ is rather simpler
than $\CS(\Psi)$ the proof is much quicker. 

\ms Our second main result concerns stochastic completeness.

\begin{definition} {\rm
The process associated $X$ is called 
{\em stochastically  complete} if $P_t 1=1$ $m$-a.e.\ for 
some (or equivalently all) $t>0$. }
\end{definition}

The energy measure $\Gamma$ defines in an intrinsic way  a 
pseudo metric  $\varrho$  on $(\sX,m)$ by
\be \label{e:intrin}
 \varrho(x,y) = \sup\{ f(y)-f(x): f\in \sF, \quad d\Gam(f,f) \le d m \} 
\ee
called the intrinsic metric or Carath\'eodory metric. 
We will denote by $B_\rho(x,r)=\{ y\in \sX:\, \varrho(x,y)<r \}$ the open 
ball with center $x$ and radius $r$ w.r.t.\ the $\varrho$ metric. 
Further, we will use the notation
\[
\varrho(x,\infty)
 :=\sup\{ r>0: B_\rho(x,r) \mbox{ is relatively compact } \subset \sX\}.
\]
If $\sX$ is a Riemannian manifold and $\sE(f,f)=\int |\grad f|^2\, d\mu$, 
then $\rho$ is just the Riemannian metric.

The pseudo-metric $\rho$ is not always useful. 
For some fractal sets such as the Sierpinski carpet
the measures $\Gam(f,f)$ and $m$ are mutually singular -- see \cite{Hi}.
In these cases the only functions $f$ satisfying the conditions of 
\eqref{e:intrin} are constants, and so $\rho$ is identically zero.

\ms
The following theorem gives, in the manifold case, the best possible 
criterion for stochastic completeness in terms of volume growth.

\begin{theorem}[\cite{Gr1, Gr2, St1}] \label{T:vgc}
Suppose that the metrics $\varrho$ and $d$ on $\sX$ are equivalent,
and all balls $B_\rho(x,r)$ are relatively compact. We say that
(VGC) holds if for some $x \in \sX$,
\begin{align} \label{e:vgc}
\int_1^\infty \frac r {\log m(B_\rho(x,r))} dr=\infty.
\end{align}
If (VGC) holds then $(\sX,\sE)$ is stochastically complete.
\end{theorem}

Our second main theorem gives a criterion for stochastic
completeness, in terms of a balance between the energy of cutoff
functions between a sequence of compact sets, and the volume of the
regions between these sets.

\begin{theorem} \label{T:ScomI}
Let $D_n$ be an increasing sequence of open sets with compact closure,
such that $\cup D_n =\sX$. Write $U_n=D_{n+1}-D_n$. Let
$\th_n>0$ be such that $\CSD(D_n, D_{n+1}, \th_n)$ holds for each $n$. \\
(a) Suppose that $\th_n \le c_1$ for all $n$. If
\be  \label{e:asybI}
 \liminf_n \frac{ \th_n m(U_n)}{4^n} =0
\ee
then stochastic completeness holds.\\
(b) Suppose $\th_n = c_0^2 n^2$, and there exists a constant $b>0$ such that
\be \label{e:asyaI}
  m(U_n) \le e^{2b (\log n)^2 }.
\ee
Then stochastic completeness holds.
\end{theorem}

\begin{remark}{\rm 
1. Note that this Theorem does not involve the intrinsic metric $\rho$. \\
2. We give an example below of a space $\sX$ such that for some
sufficiently large $R_0$ one has $\sX= B_\rho(x_0, R_0)$, but 
which is still stochastically complete. \\
3. In terms of volume growth, this Theorem gives a weaker criterion
than the results of \cite{Gr1,St1}. However since \eqref{e:asyaI}
only requires that a subsequence of annuli have small volume, there are
manifolds for which Theorem \ref{T:ScomI} gives stochastic completeness,
while the volume growth criterion of \cite{Gr1,St1} fails. \\
4. The constant 4 in \eqref{e:asybI} is not best possible; it is related to
the choice of $1/8$ in \eqref{e:csa1}. \\
5. See Remark \ref{R:scgam} for the case $\th_n \asymp n^{2\gam}$ for
$\gam \in (0,1)$.
} \end{remark}

The layout of this paper is as follows. In Section 2 we show how
$\CSA(\Psi)$ can be used to give a generalization to the space-time
scaling $\Psi$ of the `Davies-Gaffney' bound obtained by Davies
in \cite{D}. In Section 3 we use $\CSA(\Psi)$ to obtain a 
Cacciopoli type inequality. This is then used in Section 4 
to obtain mean value inequalities, which lead to 
the upper heat kernel bound $\UHK(\Psi)$. In Section 5
we prove that $\UHK(\Psi)$ implies $\CSA(\Psi)$. Section 6
proves Theorem \ref{T:ScomI}, and Section 7 gives examples, based
on the `pre-Sierpinski carpet', of 
spaces which are geodesically incomplete, or for which the criterion
of \cite{Gr1, St1} fails, but which are still stochastically complete.

\sms
We write $c$, $c'$ to denote positive constants which may change
on each appearance. Constants denoted $c_i$ will be the same through
each argument. Constants related to fundamental properties of the
space $\sX$ or Dirichlet form, such as those in the volume
doubling property, will be denoted $C_\cdot$ and will 
be the same throughout each argument.

\medskip
\noindent {\bf Acknowledgment. }
The authors wish to thank
Rich Bass for several conversations on the topic of Remark \ref{R:main}(2).

\section{ Davies Gaffney estimate}  \label{s:DG}

We begin by noting the following Cauchy-Schwarz inequality.
Let $u,v \in \sF$, $f,g\in L^\infty(\sX,m)$,  and $\lam>0$. Then
\begin{align} \nn
\int_\sX fg \, d\Gamma(u,v) 
&= \int_\sX (f/\lam^{1/2}) (\lam^{1/2}g ) \, d\Gamma(u,v) \\
\nn
&\leq \big( \lam^{-1} \int_\sX f^2 d\Gamma(u,u) \big)^{1/2} 
\cdot \big( \lam \int_\sX g^2 d\Gamma(v,v) \big)^{1/2} \\
\label{e:cs-lam}
& \le \frac{1}{2\lam} \int_\sX f^2 \, d\Gamma(u,u)
 + \frac{\lam}{2}   \int_\sX g^2 \, d\Gamma(v,v).
\end{align}

Let $D_n, n \ge 0$ be an increasing sequence of open subsets of $\sX$ with
$\ol D_n \subset D_{n+1}$. 
Suppose that 
$\CSD( D_n , D_{n+1}, \th_n)$ holds for each $n$, and let $\vp_n$
be the associated cutoff functions. 
Let $(a_n, n \ge 0)$ be an increasing sequence, with $a_0 \ge 0$.
Set
\begin{align} \label{e:vpdef0}
  \vp &= a_0 + \sum_{n=0}^\infty (a_{n+1}-a_n) (1-\vp_n), \\
  b_n &= \frac{ (a_{n+1}-a_n)}{a_n}, \qq b^* = \sup_n b_n, \\
\label{e:C0def}
  C_0 &= \sup_n b_n^2 \th_n. 
\end{align}

\begin{lemma}\label{L:CS1}
Suppose $D_n$, $\vp_n$ and $\vp$ are as above. 
Then for any  $u \in \sF$
\be \label{e:csb1}
 \int_{\sX} u^2 \, d\Gamma(\vp,\vp) 
\le \frac{(b^*)^2}{8} \int_\sX \vp^2 \, d\Gamma(u,u)
  + C_0  \int_{\sX} \vp^2 u^2 dm.
\ee
\end{lemma}

\proof
Let $U_n = D_{n+1} - D_n$, and  
note that $a_n \le \vp \le a_{n+1}$ on $U_n$.
Since $\Gam(\vp_n, \vp_m)=0$ if $n \neq m$, using CSD,
\begin{align*}
  \int u^2 \, d\Gamma(\vp,\vp) 
 &= \sum_n  (a_{n+1}-a_n)^2 \int_{U_n} u^2 \, d\Gamma(\vp_n,\vp_n)  \\
&\le  \sum_n  (a_{n+1}-a_n)^2  
\Big( {\fract18} \int_{U_n} d\Gamma(u,u) + \th_n \int_{U_n} u^2 dm   \Big) \\
&\le \sum_n  b_n^2  \Big(   {\fract18} \int_{U_n} \vp^2 \,  d\Gamma(u,u) 
  + \th_n \int_{U_n} \vp^2 u^2 dm   \Big),
\end{align*}
proving \eqref{e:csb1}. \qed

We can use this to obtain an analogue of Lemma 1 of \cite{D}.

\begin{proposition} \label{p:CSDG}
Let $\vp$, $a_n$, $b_n$, $\th_n$ and $C_0$  be as above, 
and suppose that  $b^* \le 1$. 
Let $f$ have compact support. Set $u_t = P_t f$. Then
\begin{align} \label{e:davp}
 &\| u_t \vp \|_2 \le   \| f \vp \|_2 \exp( 2 C_0 t), \\
\label{e:Iin2}
  &\int_0^t  \int \vp^2 \,  d\Gamma(u_s,u_s) ds \le 2  \| f \vp \|_2^2 e^{4C_0 t}.
\end{align}
\end{proposition}

\proof Given Lemma \ref{L:CS1} the proof is as in \cite{D}.
Let $f$ have compact support and $u_t = P_t f$.
Let $N \ge 1$, and set
\begin{align*} 
  \wt \vp_N &= a_0 + \sum_{n=0}^N (a_{n+1}-a_n) (1-\vp_n),  \\
  h_N(t) &= || u_t \wt \vp_N ||_2^2 = \int u_t^2 \wt \vp^2_N dm.
\end{align*} 
Then since $u_t \in \sD(\sL)$, $\wt \vp_N^2 u_t \in \sF$, 
\begin{align} \nn
 h_N'(t) &=  2 \langle \sL u_t, \wt \vp_N^2 u_t \rangle
 = -2\sE( u_t, \wt  \vp_N^2 u_t)\\
\label{e:hub1} 
 &=-2 \int_\sX \wt \vp_N^2 \, d\Gamma(u_t,u_t) 
 -4 \int_\sX \wt \vp_N u_t \, d\Gamma(u_t,\wt \vp_N). 
\end{align}
Using \eqref{e:cs-lam} with $\lam=2$ to bound the second term,
\be \label{e:da2}
h_N'(t) \le - \int_\sX \wt \vp_N^2\,  d\Gamma(u_t,u_t)
 + 4 \int_\sX u_t^2 \,  d\Gamma(\wt \vp_N,\wt \vp_N). 
\ee
So by Lemma \ref{L:CS1}, 
\begin{align} \nn
h_N'(t) &\le -(1- \half (b^*)^2) \int_\sX \wt \vp_N^2 \,  d\Gamma(u_t,u_t)
+ 4 C_0  \int_{\sX} \wt \vp_N^2 u_t^2 dm \\ \label{e:hprime}
&\le - \half \int \wt \vp_N^2\,  d\Gamma(u_t,u_t) +  4C_0 h_N(t).
\end{align}
Thus $h_N' \le 4C_0 h_N$, and hence
$h_N(t) \le h_N(0) \exp( 4 C_0 t)$. 
Integrating \eqref{e:hprime} we obtain
$$ h_N(t) -h_N(0) + \half \int_0^t  \int \wt \vp_N^2 \, d\Gamma(u_s,u_s) \, ds
 \le   \| f \wt \vp_N \|^2_2 \left( e^{4C_0 t} -1 \right). $$
Letting $N \to \infty$ gives \eqref{e:Iin2} and 
\eqref{e:davp}. \qed

We can use this to obtain a generalization of the `Davies-Gaffney'
bound in \cite{D}.

\begin{proposition}\label{P:DG}
Suppose $\CSA(\Psi)$ holds. Let
$x_1, x_2 \in \sX$ and let $d(x_1,x_2)=R$.
If $f_i\in L^2$ have support in $A_i= B(x_i, R/4)$ then 
\be \label{e:anDG}
 \langle P_t f_1, f_2 \rangle \le c_1 
 ||f_1||_2 ||f_2||_2 \exp( - c_2 \Phi(R,t) ). 
\ee
Here $c_i$ depend only on $C_S$, $\beta$ and $\beta_L$.
\end{proposition}

\proof 
First note that 
\be \label{e:P22}  
\langle P_t f_1, f_2 \rangle \le \| P_t f_1 \|_2 \|f_2\|_2
 \le ||f_1||_2 ||f_2||_2; 
\ee
adjusting the constants $c_i$ this is enough to give \eqref{e:anDG}
if $\Phi(R,t)$ is small.
Next, it is enough to prove \eqref{e:anDG} when $||f_i||_2=1$,
so we assume this.

Choose $m \ge 1$, 
let $r = R/2m$, and $D_k = B(x_1, \frac14 R + kr)$, $k\ge 0$.
Let $\vp_k$ be cutoff functions for $D_k \subset D_{k+1}$, for
$k \ge 0$. 
Set $a_k = 2^{ k \wedge m}$, and define $\vp$ as in \eqref{e:vpdef0}. 
Note that $b^* = (2-1)^2=1$, $\th_k = C_S \Psi(r)^{-1}$, and so
$C_0 =  C_S \Psi(r)^{-1}$.

Then writing $u_t = P_t f_1$, as in \cite[Theorem 2]{D} we have
\begin{align} \nn
  \langle P_t f_1, f_2 \rangle &=  \langle \vp P_t f_1, \vp^{-1} f_2 \rangle \\
\nn
&\le \| \vp u_t \|_2 \| \vp^{-1} f_2 \|_2 \\
\nn
&\le \| \vp f_1 \|_2 \exp(2 C_S \Psi(r)^{-1} t ) \| \vp^{-1} f_2 \|_2 \\
&\le  \exp(2 C_S  \Psi(r)^{-1} t ) ( \sup_{A_1} \vp )(\sup_{A_2} \vp^{-1}). 
\end{align}
The construction of $\vp$ gives
$\vp =1 $ on $A_1$, and $\vp = 2^m$ on $A_2$.
So
\be \label{e:Jm}
\log  \langle P_t f_1, f_2 \rangle \le -\Big(m \log 2 - \frac{ 2 C_S t}{ \Psi(R/m)} 
\Big) = -\log 2 \Big( m -  \frac{ c_3 t}{ \Psi(R/m)} \Big). 
\ee

It remains to choose $m \in \bN$ so as to obtain the bound \eqref{e:anDG},
and we need to consider several cases.

\noindent Case 1. $t\ge \Psi(R)$. By Lemma \ref{L:Phi}
we have $\Phi(R,t) \asymp 1$, and adjusting the constant $c_1$ we obtain
\eqref{e:anDG} from \eqref{e:P22}.

\noindent Case 2. $R\le  t \le \Psi(R)$. 
Then $R \ge 1$, so we have $R \le t < R^\beta$.
 We will choose $m\le R$, so the final term in \eqref{e:Jm} is
$$ m\Big( 1  - \frac{c_3 t m^{\beta-1}}{R^\beta} \Big). $$
We wish to choose $m$ so that $c_3 t m^{\beta-1}R^{-\beta} \in [1/3,2/3]$, and
this will be possible provided $R^\beta/t$ is greater than some constant 
$c_4$ (depending only on $c_3$ and $\beta$). We then have
$m \asymp (R^\beta/t)^{1/(\beta-1)}$, and hence we obtain the bound \eqref{e:anDG}. 
If $R^\beta/t< c_4$ then $\Phi(R,t)\le c_5$ and again we obtain 
\eqref{e:anDG} from \eqref{e:P22}.

\noindent Case 3. $t< R$ and $t \le \Psi(R)$. 
In this case we will choose $m >R$, so that $\Psi(R/m) = (R/m)^{\beta_L}$.
If $R\le 1$ then $\Psi(R)= R^{\beta_L}$ and 
so the argument is as in Case 2. If $R>1$ then $t<R < R^{\beta_L}$,
so again we can proceed as in Case 2. \qed

\section{Cacciopoli and mean value inequalities} \label{s:cac}

In this section we prove a mean value inequality as in \cite[Section 3]{Gr0}.
We begin by seeing that $\CSA(\Psi)$ enables us to prove \
a Cacciopoli inequality similar to \cite[Lemma 3.1]{Gr0}. 
To that aim we need to give a definition of caloric functions 
in the general context of metric measure spaces.

\begin{definition} {\rm 
Let $I$ be an interval in $\mathbb{R}$. We say that a function 
$u:I \rightarrow L^2( \sX,m)$ is weakly differentiable at $t_0\in I$ if for any 
$f\in L^2(\sX,m)$ the function $\langle u (t) , f \rangle$ is differentiable 
at $t_0$. By the principle of uniform boundedness, in this case there is a function
$w \in L^2(\sX,m)$  such that
\begin{align*}
\lim_{t\to t_0} \left( \frac{u (t)-u (t_0) }{t-t_0} , f \right) =\langle w,f\rangle
\end{align*}
for all $f\in L^2(\sX,m)$. We refer to the function $w$ as the weak derivative 
of the function $u$ at $t_0$ and write 
$w=\frac{\partial}{\partial t} u(t_0)=u_t(t_0)$.
} \end{definition}

\begin{definition} {\rm 
Consider a function $u : I \rightarrow \sF$  and let $\Omega$ be an 
open subset of $\sX$. We say that $u$ is a caloric function in $I\times \Omega$ 
if $u$ is weakly differentiable in the space $L^2(\Omega)$ at any $t\in I$ and, 
for any non-negative $f\in \sF_\Omega$ and for any $t\in I$,
\begin{align*}
\langle u_t , f \rangle + \sE (u, f ) = 0.
\end{align*}	
} \end{definition}

\begin{lemma} \label{L:ci}
(Cacciopoli inequality.) 
Let $x_0 \in \sX$, $B= B(x_0, R)$, $r<R$ and $B'=B(x_0, R-r)$.
Suppose that $\CSD(B', B, \th)$ holds, and let $\vp$ be
the associated cutoff function for $B' \subset B$.
Let $T>0$, and set $Q=B \times (0,T)$.
Let $k(t)$ be a Lipschitz function of $t$ with $k(0)=0$,
$0 \le k \le 1$ and $||k'||_\infty =K$. 
Let $u=u(x,t)$ be a non-negative caloric function, and 
$v = (u-\th)_+$, where $\th>0$. Set
$$ \eta(x,t) = \vp(x) k(t). $$
Then 
\be \label{e:cacc}
 \int_B v(x,T)^2 \eta(x,T)^2 m(dx) 
  + \frac{2}{9} \int_Q d\Gam( \eta v , \eta v) \, dt
\le 2( \frac{20} 9 \th  + K) \int_Q v^2 \, dm \, dt .
\ee
\end{lemma}

\proof 
Since $k(0)=0$ we have, writing $v_t = \pd v/\pd t$, 
\be \label{e:cac1}
  \frac12 \int_B v(x,T)^2 \eta(x,T)^2 m(dx) =  \int_Q v v_t \eta^2 dm\, dt +
  \int_Q v^2 \eta \eta_t \, dm \, dt. 
\ee
Using the fact that $u$ is caloric we get
\begin{align} \label{e:cac2}
 \int_Q \eta^2 v v_t \, dm\, dt 
  &=\int_Q \indicator_{\{u>\th\}} \eta^2 v u_t \, dm\, dt
      =  - \int_Q d\Gamma(\eta^2 v,u) \nonumber \\
 &=  - \int_Q \eta^2 \, d\Gam(v,v)\, dt  - 2 \int_Q v \eta \, d\Gam(v, \eta)\, dt . 
\end{align}
Further,
\be \label{e:cac3}
 \int_B d\Gam(v \vp, v\vp) \le
2\int_B \vp^2 \, d\Gam(v,v) + 2 \int_B v^2 \,  d\Gam(\vp,\vp).
\ee
Let $\lam>0$. Then using \eqref{e:cs-lam}
\begin{align*} 
  - \int_B \vp^2 \, d\Gam(v,v)  - 2 \int_B v \vp \, d\Gam(v, \vp)
\le  (-1 + \lam^{-1})  \int_B \vp^2 d\Gam(v,v) + \lam  \int_B v^2 d\Gam(\vp, \vp).
\end{align*}
Taking $\lam=2$ and using \eqref{e:cac3} and CSD we obtain
\begin{align} \nn 
  - \int_B \vp^2 \, d\Gam(v,v)  - &2 \int_B v \vp \, d\Gam(v, \vp) 
  + a \int_B d\Gam(v \vp, v\vp)  \\
\nn
&\le (-\frac12 + 2a) \int_B \vp^2 \, d\Gam(v,v) + (2+ 2a) \int_B v^2 d\Gam(\vp, \vp)\\
\nn
&\le  (-\frac12 + 2a + (2+2a) \frac18 )  \int_B \vp^2 d\Gam(v,v) +
   (2+2a) \th \int_B v^2  \, dm\\
\label{e:fr18}
&= \frac {20}{9} \th \int_B v^2 \, dm
\end{align}
if $a=1/9$. Multiplying this inequality by $k(t)^2$ and integrating gives
\begin{align*} 
   - \int_Q \eta^2 \, d\Gam(v,v)\, dt  &- 2 \int_Q v \eta \, d\Gam(v, \eta)\, dt 
  + \frac19 \int_Q d\Gam(v \eta, v\eta) \,dt  \\
&\le \frac {20} 9 \th   \int_0^T \int_B v^2 k(t)^2 \,dm\, dt 
\le \frac {20} 9 \th \int_Q v^2 \, dm\, dt .
\end{align*}
Combining this with \eqref{e:cac1} and \eqref{e:cac2} we obtain
\begin{align*}
  \int_B v(x,T)^2 \eta(x,T)^2 \, m(dx) + \frac29  &\int_Q d\Gam(v \eta, v\eta)\, dt  \\
&\le 2 \int_Q v^2 \eta \eta_t \, dm \, dt 
           +  \frac{40} 9 \th \int_Q v^2 \, dm \, dt \\
 &\le  (2K + \frac{40} 9 \th ) \int_Q v^2 \, dm \, dt .
\end{align*}
In the final line we used the fact that $\eta\le 1$. \qed

\begin{remark} {\rm 
Note that to obtain \eqref{e:fr18} we needed that the constant in the 
first term on the right of \eqref{e:csa1} was less than $1/4$. 
} \end{remark}

The key step in the proof of the mean value inequality is the following
comparison over cylinders. For a cylinder $Q \subset \sX \times \bR_+ $
and a function $w$ write
$$ I(w,Q) = \int_Q w^2 \, dm \, dt. $$

\begin{lemma} (See \cite[Lemma 3.2]{Gr0}.)
Suppose $\FK(\Psi)$ and $\CSA(\Psi)$ hold. 
Let $\wt u$ be a caloric function in $Q =  B(x_0,R) \times (0,T)$. Let $u=\wt u_+$,
$\th > 0$ and
$v = (u-\th)_+$. Let $0<T_1<T$, $R_1\in (\half R,R)$, 
$Q_1 = B(x_0,R_1)\times (T_1,T)$, 
$$ I = I(u,Q), \q I_1=I(v,Q_1), $$
and $\delta= T_1 \wedge \Psi(R-R_1)$. Then
\be \label{e:IIi}
 I_1\le 
\frac{ c_1 I^{1+\nu} \Psi(R) }{ \delta^{1+\nu} \th^{2\nu} \left. m(B)\right.^\nu}.
\ee
\end{lemma}

\proof
Let $B=B(x_0,R)$, $B' = B(x_0, \half(R+R_1))$ and $B_1=B(x_0,R_1)$.
Set
$$ D_t =\{ x \in B': u(x,t) > \th \}. $$
Let $\vp$ be a cutoff function for $B_1 \subset B'$, $k(t) = 1 \wedge (t/T_1)$,
and $\eta(x,t)= \vp(x) k(t)$. 

As in \cite{Gr0} the proof uses five inequalities:
\begin{align}
\label{e:in15}
 \int_{B'} u(x,t_0)^2\,  dm 
&\le c_0 \delta^{-1} \int_Q u^2 \, dm\, dt \, \hbox{ for } t_0 \in (T_1,T),\\
\label{e:in25}
\int_Q d\Gam( v\eta, v \eta)\, dt &\le c_0 \delta^{-1} \int_Q u^2 \, dm \, dt , \\ 
\label{e:in35}
 \int_{B} d\Gam( v\eta, v \eta) &\ge \lam_1(D_t) \int_B v^2 \eta^2\, dm 
   \, \hbox{ for } t \in (0,T),\\
\label{e:in55}
 \lam_1(D_t) &\ge C_F m(B)^\nu \Psi(R)^{-1} m(D_t)^{-\nu}, \\
\label{e:in45}
 m(D_t) &\le \th^{-2} \int_{B'} u(x,t)^2 \, dm. 
\end{align}
Of these, \eqref{e:in35} is immediate from the variational definition of $\lam_1$, 
\eqref{e:in55} is the Faber-Krahn inequality \eqref{e:fki},
and \eqref{e:in45} is just Markov's inequality.
So it remains to prove  \eqref{e:in15} and  \eqref{e:in25}.

The inequality \eqref{e:in25} is immediate from \eqref{e:cacc}. 
Since $\|k'\|_\infty=1/T_1$ we have the constant on the right side of
\eqref{e:cacc} is $c( \Psi(R-R_1)^{-1} + T_1^{-1}) \le c' \delta^{-1}$. So 
$$ \frac29 \int_Q d\Gam( v\eta, v \eta)\, dt 
 \le c \delta^{-1} \int_Q u^2\, dm \,dt . $$
For  \eqref{e:in15} let $\wt \vp$ be a cutoff function for $B' \subset B$
and $\wt \eta (x,t)= \wt \vp(x) k(t)$. Then by  \eqref{e:cacc} applied to $u$
in the cylinder $Q_t= B\times (0,t)$,
$$   \int_B u(x,t)^2 \wt \eta(x,t)^2 m(dx) 
 \le c \delta^{-1}  \int_{Q_t} u^2 \, dm \, dt
 \le c \delta^{-1}  \int_{Q} u^2 \, dm \, dt. $$

The rest of the argument is as in \cite{Gr0}. 
\qed

\section{ Heat kernel upper bounds} \label{s:hku}

These bounds can now be proved by the methods of \cite{CG}, which in turn
uses ideas in \cite{Gr0}.
Since \cite{CG} is written in the graph context, and both of these
papers just consider the case $\Psi(r)=r^2$, we give details.
In particular we need to be more careful in our handling of 
exceptional sets; issues with these do not arise for
the manifolds or graphs treated in \cite{Gr0, CG}.
Note that VD implies that there exists a constant $\al<\infty$ such that
\be \label{e:vdtr}
 \frac{V(x,R)}{V(y,r)} \le C_G \Bigl( \frac{d(x,y)+R}{r} \Bigr)^{\al},
\q 0<r<R, \, x,y \in \sX.
\ee
Define the measure $\wt m(dx,ds) = m(dx) ds$ on $\sX \times \bR$.
Given a cylinder $Q \subset \sX \times \bR$ and $u : Q \to \bR$
we write $\esssup {u}$
for the essential supremum with respect to the measure $\wt m$.

Define 
\be \label{e:b1b2}
\beta_2 = \beta_L \vee \beta, \q \beta_1 = \beta_L \wedge \beta,
\ee 
and note that if $0<r<R$,
\be \label{e:Psig}
 \Big( \frac{R}{r} \Big)^{\beta_1} \le \frac{\Psi(R)}{\Psi(r)}
  \le \Big( \frac{R}{r} \Big)^{\beta_2}, \q
 \Big( \frac{R}{r} \Big)^{1/\beta_2} \le \frac{\Psiinv(R)}{\Psiinv(r)}
  \le \Big( \frac{R}{r} \Big)^{1/\beta_1}.
\ee
Write
\be
 F(R,T) = \frac{T}{\Psi(R)} \vee \Big( \frac{\Psi(R)}{T} \Big)^{\al/\beta_1}. 
\ee

\begin{proposition} \label{P:mvi2}
($L^2$ Mean value inequality).
Set $Q=B(x_0,R) \times (0,T)$.
Assume $\CSA(\Psi)$ and $\FK(\Psi)$ hold, and let $u \ge 0$ 
be caloric in $Q$. Then if $Q_\infty = B(x_0,R/2) \times (T/2,T)$,
\be \label{e:mvi2}
 \esssup_{Q_\infty} u(x,s)^2  \le  \frac{c_0 F(R,T) }{T V(x_0,R)} 
    \int_Q u(x,s)^2 \wt m(dx,ds).
\ee
\end{proposition}

\proof (See the proof of \cite[Theorem 3.1]{Gr0}.)
It is sufficient to consider the case $T=\Psi(R)$. Indeed, suppose \eqref{e:mvi2}
holds in this case, and let $T=\lam \Psi(R)$.

If $\lam \in (0,1)$ let $r$ be such that $\Psi(r) =T$.
We can cover $B'=B(x_0, R/2)$ by balls $B(z_i,r/2)$ such that 
each $B(z_i,r) \subset B(x_0,R)$. Let $Q_i=B(z_i,r) \times (0,T)$,
and $Q_{i, \infty}=B(z_i,r/2) \times (T/2,T)$.
Note that by \eqref{e:Psig} 
\be \label{e:vdtr2}
 \frac{V(x_0,R)}{V(x_0,r)} \le C_G \Big( \frac{R}{r}\Big)^\al 
 = C_G \Big( \frac{ \Psiinv(T/\lam)}{\Psiinv(T)}\Big)^\al 
\le C_G \lam^{-\al/\beta_1}.
\ee

Then
$$  \esssup_{Q_\infty} u(x,s)^2  \le \max_i  \esssup_{Q_{i,\infty}} u(x,s)^2, $$ 
and for each $i$, using \eqref{e:vdtr2},
\begin{align*} 
 \esssup_{Q_{i,\infty}} u(x,s)^2 
 &\le \frac{c_0}{T V(x_0,r)} \int_{Q_i} u^2 d\wt m  \\
 &\le  \frac{c_0}{T V(x_0,R)} \frac{V(x_0,R)}{V(x_0,r)} \int_{Q} u^2 d\wt m \\
 &\le  \frac{c_0C_G  \lam^{-\al/ \beta_1}}{T V(x_0,R)} \int_{Q} u^2 d\wt m.
\end{align*}
Similarly if $T=\lam \Psi(R)$ with $\lam>1$ then applying \eqref{e:mvi2}
to a sequence of cylinders $(T_i-\Psi(R),T_i)\times B(x_0,R)$ one obtains
$$  \esssup_{Q_\infty} u(x,s)^2 
\le  \frac{c_0 \lam }{T V(x_0,R)} \int_Q u(x,s)^2 d\wt m. $$

Now let $T=\Psi(R)$. Let $\delta_k$, $k=0, 1, \dots$ be a sequence to
be chosen later, and let $(r_k)$ and $(t_k)$ be sequences such that
$r_0=R$, $t_0=0$,
\be
 \delta_{k+1} = \Psi(r_k- r_{k+1})= t_{k+1} - t_{k},
\ee
and 
\be\label{e:rtcond}
R =r_0 > r_1 > \dots > r_k > \dots > R/2, \q  0=t_0 < t_1 < \dots < t_k< \dots < T/2.
\ee

Let $Q_k = B(x_0,r_k) \times (t_k,T)$.
Let $\th>0$ (also to be chosen later) and set
$\al_k = (1- 2^{-k})\th$, $u_k = (u- \al_k)_+$, and
$$ I_k = \int_{Q_k} u_k^2 d\wt m. $$
Let 
$\th_k = \al_{k+1}-\al_k$. Then by \eqref{e:IIi} applied to the function
$u_k$ in $Q_{k+1} \subset Q_k$,
\be
 I_{k+1} \le I_k \;
 \frac{ c_1 I_k^\nu \Psi(r_k) }{ \delta_{k+1}^{1+\nu} \th_k^{2\nu} m(B(x_0,r_k))^\nu}
\le I_k \; 
 \frac{ c_2 \Psi(R)  2^{2\nu(k+1)} I_k^\nu } 
   {V(x_0,R)^\nu \th^{2\nu} \delta_{k+1}^{1+\nu} }.
\ee
Let $A= c_2 \Psi(R)  4^\nu/ V(x_0,R)^\nu  \th^{2\nu} $.
Now choose $M\ge 2$ and $b>4$ 
such that $M^{1/\beta_2} \ge 2$, and  $(4/b)^{1/(1+\nu)} = M^{-1}$. 
Choose $\delta_{k+1}$ such that 
\be
 \frac{A 2^{2\nu k }I_k^\nu}{\delta_{k+1}^{1+ \nu}} = b^{-1}, \q k=0, 1, \dots 
\ee
With this choice of $(\delta_k)$ we have $I_k \le b^{-k} I_0$, and consequently
\be \label{e:dkb}
 \delta_{k+1} = (A b 4^{k\nu} I_k^\nu)^{1/(1+\nu)} \le (A b I_0^\nu)^{1/(1+\nu)} 
  (4/ b)^{k \nu/(1+\nu)} =  A_1 M^{-k} ,
\ee
where $A_1 =  (A b I_0^\nu)^{1/(1+\nu)}$.
In order that the condition \eqref{e:rtcond} should hold, we need 
\begin{align}\label{e:tc1}
 &\sum_{k=1}^\infty \delta_{k} \le T/2, \\
\label{e:rc1}
 &\sum_{k=1}^\infty \Psiinv(\delta_{k}) \le R/2.
\end{align}
We have
$$ \sum_{k=1}^\infty \Psiinv(\delta_{k}) \le \sum_{k=1}^\infty \Psiinv(A_1 M^{-k})
 \le c_3 \Psiinv(A_1) \sum_{k=1}^\infty M^{-k/\beta_2} \le 2c_3  \Psiinv(A_1), $$
and 
$$ \sum_{k=1}^\infty \delta_{k} \le 2 A_1.  $$
So for \eqref{e:rtcond} to hold it is enough that
$$ A_1=  \Big( \frac{c_24^\nu b I_0^\nu \Psi(R)}{ V(x_0,R)^\nu \th^{2 \nu}}\Big)^{1/(1+\nu)}
  \le (T/4) \wedge \Psi(R/4c_3). $$
For this it is enough if $\th$ is chosen large enough so that
$$ \frac{ I_0^\nu \Psi(R)}{ V(x_0,R)^\nu \th^{2 \nu}} \le c_4 \Psi(R)^{1+\nu}, $$
and so we can take
\be 
\th^2 = \frac{c_5 I_0}{\Psi(R) V(x_0,R)}.
\ee
We then have
$I_k \to 0$ as $k \to \infty$, and hence 
\be
 \int_{Q_\infty} (u(x,s)- \th)_+^2 \wt m(dx,ds) \le \inf_k I_k =0,
\ee
which implies that $u(x,s) \le \th$ $\wt m$ a.e. on  $Q_\infty$. \qed

\sms
We now give an $L^1$ mean value inequality.

\begin{proposition} \label{P:mvi}
($L^1$ Mean value inequality).
Assume $\FK(\Psi)$ and $\CSA(\Psi)$ hold. Let $R>0$, $T= \Psi(R)$, let
$Q=B(x_0,R) \times (0,T)$, and let $u \ge 0$ be caloric in $Q$. Then
writing $Q' = B(x_0,2R/3) \times (T/2,T)$,
\be \label{e:mvi}
 \esssup_{Q'} u(x,s)  \le  \frac{c_1 F(R,T) }{T V(x_0,R)} 
    \int_Q u(x,s) \wt m(dx,ds).
\ee
\end{proposition}

\proof 
This follows from the $L^2$ mean value inequality by quite general arguments,
which use only VD --  see p.\ 688-691 of \cite{CG}. 
As with Proposition \ref{P:mvi2}, it is enough to consider the case $T=\Psi(R)$.
\qed

In order to obtain heat kernel bounds from the mean value theorem, we
need better control of the exceptional set. We will use regularity
results from \cite{GT}, and to use these we need to consider the 
killed heat kernel. For $D \subset \sX$ write  
$(P^{D}_t)$ for the semigroup of $X$ killed on exiting $D$. Then if
$\FK(\Psi)$ holds, by 
\cite[Lemma 5.5]{GH} the semigroup  $(P^{D}_t)$ is ultracontractive,
that is there exists a left continuous function $\gamma(t)$ such that
$$ || P^D_t f ||_\infty \le \gam(t) ||f||_1, \q t >0, f \in L^1(\sX)\cap L^2(\sX). $$
(In fact we have $\gam(t) = c (at)^{-1/\nu}$ with $a=a(n) = m(D_n)^\nu/\Psi(nR)$).
Consequently we will be able to use \cite[Theorem 2.12]{GT} to obtain 
estimates which hold on $\sX-\sN$, where $\sN$ is a properly exceptional set.

\begin{lemma} \label{L:mvpT}
Assume $\FK(\Psi)$ and $\CSA(\Psi)$. 
Let $x_0, x_1 \in \sX_0$, 
$T>0$, $0<t \le T$, and $r = \Psi(t)$. 
Let $R_0 > R$ and $D=B(x_0, R_0)$. Then for $m$-a.a. $y \in B(x_0,r/2)$,
\be \label{e:mvpT}
 p^D_T(x_1,y) \le \frac{c_1}{tV(x_0,r)} 
 \int_{T-t/2}^{T+t/2} \int_{B(x,r)} p^D_s(x_1,y') \wt m(dy',ds).
 \ee
\end{lemma}

\proof 
Set $Q=B(x_0,r) \times (T-t/2,T+t/2)$, and $Q'=B(x_0, r/2) \times (T,T+t/2)$.
Since $u(y,s)= p^D_s(x_1,y)$ is caloric in $Q$, by Proposition \ref{P:mvi}
\be \label{e:lreg1}
 \esssup_{(y,s) \in Q'} p^D_s(x_1,y) \le
 \frac{c_1}{t V(x_0,r)}\int_Q p^D_s(x_1,y') \wt m(dy',ds).
\ee
Setting
\be
 A=  \frac{c_1}{t V(x_0,r)}\int_Q p^D_s(x_1,y') \wt m(dy',ds),
\ee
the right side of \eqref{e:lreg1} is bounded by $A$.

Thus there exists a subset $I \subset (T, T+t/2)$ of full measure
such that if $T'\in I$ then
$$ p^D_{T'}(x_1, y) \le A, \q \hbox{ for $m$-a.a. $y \in B(x_0, r/2)$ }. $$ 
Write $g(y)=p^D_T(x_1, y)$. Then the $L^2$-continuity of $(P^D_t)$ implies
that $ P^D_h g \to g$ $m$-a.e. as $h \to 0$. Taking the limit along a
sequence $h_k$ such that $T+h_k \in I$ for each $k$, it follows that
$p^D_{T}(x_1, y) \le A$ for $m$-a.a. $y \in B(x_0, r/2)$.
\qed

\begin{theorem} \label{T:hkub}
Assume VD,  $\tFK(\Psi)$, and $\CSA(\Psi)$ hold.
Then UHK($\Psi$) holds.
\end{theorem}

\proof We use the argument of \cite{CG}, but need extra care because
of exceptional sets. 
Fix $x_0,y_0 \in \sX_0$, $T>0$, and let $R= d(x,y)$. 
Let $R_0 > 4R$, and let $D= B(x_0, R_0)$.
Set $T_0=T/2$ and $r=\Psiinv(T_0)$. 
Let $Q(z)= B(z,r) \times (T-T_0/2, T+T_0/2)$. 
Let $r' < R/4 \wedge r/2$, and 
$g_1$ and $g_2$ be non-negative bounded functions with supports
in $B(x_0, r')$ and $B(y_0, r')$ respectively, such that
$\int g_1 = \int g_2 =1$. Set
$$ I = \iint p^D_T(x,y) g_1(x) g_2(y) m(dx) m(dy). $$

Let $x \in B(x_0, R/4)$. Then applying \eqref{e:mvpT} to the caloric 
function $u(y,s) = p^D_s(x,y)$ in $Q(y_0)$, we have 
\be \label{e:mv-a}
   p^D_T(x,y)  \le  \frac{c_1}{T_0 V(y_0,r)}
  \int_{T-T_0/2}^{T+T_0/2} \int_{B(y_0,r)} p^D_s(x,y') \wt m(dy',ds)
\ee
for $m$-a.a. $y \in B(y_0, r/2)$.
Hence
\begin{align}\nn
 I &\le   \frac{c_1}{T_0 V(y_0,r)}  \iint g_1(x) g_2(y)m(dx) m(dy)
       \int_{Q(y_0)} p^D_s(x,y') \wt m(dy',ds)  \\
\label{e:mv-c}
 &=  \frac{c_1}{T_0 V(y_0,r)} 
       \int_{B(y_0,r)} \int_{T-T_0/2}^{T+T_0/2} \int g_1(x) p^D_s(x,y') \wt m(dy',ds)m(dx).
 \end{align}
If $y' \in B(y_0, r)$ and $s \in (T-T_0/2, T+T_0/2)$, then 
by considering the cylinder $B(x_0,r)\times (s-T_0/2, s+T_0/2)$, we have
by \eqref{e:mvpT}, for $m$-a.a. $x \in B(x_0, r/2)$,
\begin{align*} 
   p^D_s(y',x)  &\le  \frac{c_1}{T_0 V(x_0,r)}
  \int_{s-T_0/2}^{s+T_0/2} \int_{B(x_0,r)} p^D_{s'}(y',x') \wt m(dx',ds') \\
  &\le  \frac{c_1}{T_0 V(x_0,r)}
  \int_{T-T_0}^{T+T_0} \int_{B(x_0,r)} p^D_{s'}(x',y') \wt m(dx',ds'),
\end{align*}
Substituting this into the final term in \eqref{e:mv-c}, we obtain
\begin{align} \nn
 I &\le  \frac{c_1^2}{T_0^2 V(y_0,r)V(x_0,r)} 
       \int_{B(y_0,r)} m(dy') \int_{T-T_0/2}^{T+T_0/2}ds  \\
\nn
 & \qq \qq \times \int m(dx) g_1(x) 
  \int_{T-2T_0}^{T+2T_0} ds' \int_{B(x_0,r)}m(dx') p^D_{s'}(x',y') \\
\nn
 &=  \frac{c_1^2}{T_0 V(y_0,r)V(x_0,r)} 
      \int_{B(y_0,r)} m(dy') \int_{B(x_0,r)}m(dx') 
    \int_{T-T_0}^{T+T_0} ds' p^D_{s'}(x',y') \\
\label{e:pt3i}
  &\le  \frac{2 c_1^2}{T V(y_0,r)V(x_0,r)}  \int_{T/2}^{3T/2}
  \langle P_s f_x, f_y \rangle ds,
\end{align} 
where $f_x = 1_{B(x_0,r)}$ and $f_y = 1_{B(y_0,r)}$.

If $r< R/4$ then
the Davies Gaffney bound Proposition \ref{P:DG} implies that
for $T/2 \le s \le 2T$,
\begin{align}\nn
\langle P_s f_x ,f_y\rangle 
 &\le c_2 V(x_0,r)^{1/2} V(y_0,r)^{1/2} \exp( -c_2 \Phi(R,s)  ) \\
\label{e:dga1}
 &\le c_2 V(x_0,r)^{1/2} V(y_0,r)^{1/2} \exp( -c_3 \Phi(R,T)).
\end{align}
If $r\ge R/4$ we still have by \eqref{e:P22}
\be 
 \langle P_s f_x ,f_y\rangle \le ||P_t f_x||_2 ||f_y||_2 
\le V(x_0,r)^{1/2} V(y_0,r)^{1/2}. 
\ee
As $r\ge R/4$ we have $T > c \Psi(R)$ and so 
$\Phi(R,T) \le c'$ by Lemma \ref{L:Phi}, and 
the exponential in \eqref{e:dga1} is of order 1. 
Adjusting the constant $c_2$ we therefore obtain, in both cases,
\be \label{e:hk-4}
 \langle P_s f_x ,f_y\rangle   
\le  c_2 V(x_0,r)^{1/2} V(y_0,r)^{1/2} \exp( -c_3 \Phi(R,T)).
\ee

Now set
$$ \wt V(x,r) = r^{-1} \int_r^{2r} V(x,s) ds. $$
Then $\wt V(x,r) \asymp V(x,r)$, and the function $x \to \wt V(x,r)$ is
continuous. Set
$$ H_t(x,y) = 
 c_4 \wt V(x,\Psiinv(t))^{1/2} \wt V(y,\Psiinv(t))^{1/2} \exp( -c_3 \Phi(d(x,y),t)). $$
Then from \eqref{e:hk-4} and \eqref{e:pt3i} we deduce that
$ p_T^D(x,y) \le H_T(x_0,y_0)$, for $m \times m$ a.a. $(x,y)$ in a neighbourhood
of $(x_0, y_0)$. It follows that
$ p^D_T(x,y) \le H_T(x,y)$ for $m \times m$ a.a. $(x,y)$.

Now let $D \uparrow \sX$; then since $P^D$ is ultracontractive
it follows from \cite[Theorem 2.12(c)]{GT} that 
$$ p_t(x,y) \le H_t(x,y)\q \hbox{ for $m \times m$ a.a. $(x,y)$ and all $t>0$}.$$
Since the function $H_t$ is continuous, by \cite[Theorem 2.12(d)]{GT} there exists
a properly exceptional set $\sN_1$ such that if $\sX_1 = \sX-\sN_1$ then
\be
 p_t(x,y) \le H_t(x,y)\q \hbox{ for $(x,y) \in \sX_1 \times \sX_1$ 
and all $t>0$},
\ee
which proves $\UHK(\Psi)$. \qed

\section{Proof of CSA from upper heat kernel bounds} \label{s:csa}

In this section we prove the implication $(2) \Rightarrow (1)$
of Theorem \ref{T:main}. We assume throughout this section
that  $\sX$ is unbounded, and satisfies VD and $\UHK(\Psi)$.

We begin by seeing that it is enough to prove \eqref{e:csa1}
in a slightly weaker form. 

\begin{lemma} \label{L:cswk}
Let $\sX$ satisfy VD.
Suppose that there exists constants $c_1, c_2$ such that for
all $x \in \sX$, $R>0$ and $r>0$ there exists a cutoff function
$\vp$ for $B(x,R) \subset B(x,r+R)$ such that, if $U=B(x,R+r)-B(x,R)$
and $f \in \sF$, then 
\be \label{e:csa-C}
 \int_U f^2 \, d\Gamma(\vp,\vp) 
\le c_1  \int_U d\Gamma(f,f) + c_2  \Psi(r)^{-1} \int_U f^2 dm .
\ee 
Then $\sX$ satisfies $\CSA(\Psi)$. 
\end{lemma} 

\proof
Let  $x \in \sX$, $R>0$ and $r>0$ and $B'=B(x,R)$, $B=B(x,R+r)$; 
we will construct a cutoff
function $\vp$ for  $B' \subset B$ which satisfies
\eqref{e:csa1} with $\th= c \Psi(r)^{-1}$.

Let $\lam>0$, and let 
$$ s_n = c_0 r e^{- n \lam/\beta_2}, $$
where $c_0=c_0(\lam)$ is chosen so that $\sum_{n=1}^\infty s_n =r $ and $\beta_2$
is as in \eqref{e:b1b2}. 
Set $r_0=0$, 
$$ r_n = \sum_{k=1}^n s_k, $$
so that $R < R+r_1 < R+r_2 < \dots < R+r$.
Let $B_n= B(x_0, R+r_n)$, and $U_n= B_{n+1}-B_{n}$.
By hypothesis there exists a cutoff function $\vp_n$ for $B_n \subset B_{n+1}$
satisfying 
\be \label{e:cswk}
 \int_{U_n} f^2 \, d\Gamma(\vp_n,\vp_n) 
\le c_1 \int_{U_n} d\Gamma(f,f) 
  + c_2 \Psi(s_{n+1})^{-1} \int_{U_n}  f^2 dm.  
\ee

Let $b_n = e^{-n \lam}$ and set
\be
 \vp = \sum_{n=1}^\infty (b_{n-1} - b_n) \vp_n.
\ee
Then $\vp=0$ on $B^c$, and $\vp=1$ on $B'$, so $\vp$ is a cutoff
function for $B' \subset B$. On $U_n$ we have
$$ \vp = (b_{n-1}-b_n) \vp_n + b_n, $$
and so $b_n \le \vp \le b_{n-1}$ on $U_n$. Hence on $U_n$
\be \label{e:vpne}
 b_{n-1}-b_n  \le  \frac{ \vp (b_{n-1}-b_n)}{b_n} = (e^\lam -1)\vp. 
\ee
Now if $f : B \to \bR$ then by \eqref{e:cswk}
\begin{align*} 
  \int_{B} f^2 \, &d\Gamma(\vp,\vp) 
 = \sum_{k=1}^\infty (b_{k-1}-b_k)^2 \int_{U_k} f^2 d\Gam(\vp_k,\vp_k) \\
&\le c_1 \sum_{k=1}^\infty (b_{k-1}-b_k)^2 \int_{U_k}  d\Gam(f,f) 
   + c_2 \sum_{k=1}^\infty (b_{k-1}-b_k)^2 \Psi(s_{k+1})^{-1} \int_{U_k} f^2  dm .
\end{align*}
Using \eqref{e:vpne} we have
\begin{align}\nn
 c_1 \sum_{k=1}^\infty (b_{k-1}-b_k)^2 \int_{U_k}  d\Gam(f,f) 
&\le c_1 (e^\lam-1)^2 \sum_{k=1}^\infty  \int_{U_k} \vp^2 d\Gam(f,f) \\
\label{e:comb1}
&\le c_1 (e^\lam-1)^2  \int_U \vp^2 d\Gam(f,f). 
\end{align}
Now using \eqref{e:Psig} and \eqref{e:vpne}
$$ \frac{ \Psi(r)}{\Psi( s_{k+1})} 
 \le \Big(\frac{r}{c_0(\lam) r e^{- (k+1) \lam /\beta_2}}\Big)^{\beta_2}
 = \frac{ e^{\lam} e^{k\lam} } {c_0(\lam)^{\beta_2}} 
= \frac{ { e^\lam }(e^{\lam}-1) }{ c_0(\lam)^{\beta_2}(b_{k-1}-b_k )}.
$$
Therefore 
$$  (b_{k-1}-b_k) \Psi(s_{k+1})^{-1} \le  c_3(\lam) \Psi(r)^{-1}, $$
and hence 
\begin{align}\nn
 c_2 \sum_{k=1}^\infty  \int_{U_k} (b_{k-1}-b_k)^2 \Psi(s_{k+1})^{-1} f^2  dm  
&\le c_2 c_3(\lam) \Psi(r)^{-1}  
\sum_{k=1}^\infty \int_{U_k} f^2 (b_{k-1}-b_{k}) dm \\
\label{e:comb2}
 &\le c_2 c_3(\lam) \Psi(r)^{-1} \int_U f^2  (e^\lam-1) \vp dm.
\end{align}
Thus
\be
  \int_U f^2 \, d\Gamma(\vp,\vp) \le
 c_1 (e^\lam-1)^2  \int_U \vp^2 d\Gam(f,f) 
   +  c_4(\lam) \Psi(r)^{-1} \int_U f^2 \vp dm.
\ee
We now choose $\lam$ so that $c_1^2 (e^\lam-1)^2 =1/8$ and since
$\vp\le 1$ we obtain \eqref{e:csa1}. \qed

\begin{corollary} \label{C:csastab}
Let $\sX$ satisfy VD. Then the condition $\CSA(\Psi)$ is stable.
\end{corollary}

\proof Let $(\sE_i, \sF)$, $i=1,2$ be two Dirichlet forms on
$L^2(\sX, m)$ satisfying the hypothesis of Lemma \ref{L:sgam},
and suppose that $\CSA(\Psi)$ holds for $\sE_1$. Let
$B'=B(x,R) \subset B=B(x, R+r)$, and let $\vp$ be a cutoff function
for $B' \subset B$. Then by Lemma \ref{L:sgam}, if $f \in \sF$,
$U=B-B'$,
\begin{align*}
  \int_U f^2 \, d\Gamma^{(2)}(\vp,\vp) &\le C  \int_U f^2 \, d\Gamma^{(1)}(\vp,\vp) \\
&\le (C/8) \int_U \vp^2 d\Gamma^{(1)}(f,f) + C C_S \Psi(r)^{-1} \int_U f^2 dm \\
&\le (C^2/8) \int_U d\Gamma^{(2)}(f,f) + C C_S \Psi(r)^{-1} \int_U f^2 dm .
\end{align*}
Thus $(\sX, \sE_2)$ satisfies the condition \eqref{e:csa-C}
and so by Lemma \ref{L:cswk} $\CSA(\Psi)$ holds for $\sE_2$. \qed

Now let $(X_t, t\in \bR_+, \bP^x, x\in \sX)$ be the Hunt process associated
with the semigroup $P_t$ and Dirichlet form $\sE$. 
Recall the definition of $\sX_0$ from Section 1. For a set
$D \subset \sX$ define the exit time
\be
 \tau_D = \inf\{ t >0 : X_t \in D^c \}. 
\ee

\begin{lemma}\label{L:taub1}
Suppose $\sX$ satisfies VD and $\UHK(\Psi)$. 
There exists a constant $\eps>0$ such that for all $x \in \sX_0$ 
and $r>0$,
$$ \bP^x(  \tau_{B(x,r)} \le \eps \Psi(r) ) \le \eps. $$
\end{lemma}

\proof
In the case $\Psi(r)= r^\beta$ this property is denoted $P_\beta$
in \cite{GH}, and the result follows by \cite[Theorem 2.2]{GH}.
The general case is similar. \qed
 
For $D \subset \sX$, $\lam>0$ set 
$$ G^D_\lam f(x) = \bE^x \int_0^{\tau_D} e^{-\lam t} f(X_t) dt. $$

\begin{lemma} \label{L:hprop}
Suppose $\sX$ satisfies VD and $\UHK(\Psi)$. 
Let $x_0 \in \sX$, $r>0$, $R>0$, and define the annuli
$D_0 = B(x_0,R+ 9r/10 )-\ol B(x_0, R + r/10)$,
$D_1 = B(x_0,R+ 4r/5 )-\ol B(x_0, R + r/5)$,
$D_2 = B(x_0,R+ 3r/5 )-\ol B(x_0, R + 2r/5)$. 
Let $\lam = \Psi(r)^{-1}$, and set
\be
h = G^{D_0}_\lam 1_{D_1}.
\ee
Then $h$ has support $\ol D_0$, $h \in \sF_{D_0}$ and satisfies
\begin{align} 
 h(x) &\le \Psi(r)  \q \hbox{ for all $x \in \sX$}, \\
\label{e:GDlb}
 h(x) &\ge c_1^{-1}\Psi(r) \q \hbox{ for } x \in D_2\cap \sX_0 .
\end{align}
\end{lemma}

\proof
That $h \in \sF_{D_0}$ follows by \cite[Theorem 4.4.1]{FOT}.
The definition of $h$ implies that $h(x)=0$ for $x \not\in \ol D_0$,
and the upper bound on $h$ is elementary, since $h \le G^\sX_\lam 1 = \lam^{-1}$.

Now let $\eps>0$ be as in Lemma \ref{L:taub1}.
Let $r_0=r/5$, $x \in D_2$, and $B_1 =B(x,r_0) \subset D_1$.
Let $s = \eps \Psi(r_0)$, and 
$\xi_\lam$ be an exponential r.v. independent of $X$
with mean $\lam^{-1}$. 
Then
\begin{align*}
  h(x) &\geq  \bE^x \int_0^{ \xi_\lam \wedge \tau_{D_0}} 1_{D_1}(X_t)\, dt 
 \ge  \bE^x \int_0^{ \xi_\lam \wedge \tau_{B_1}} 1_{B_1}(X_t)\, dt \\
 &\ge s \bP^x(  \xi_\lam \wedge \tau_{B_1} \ge s) 
 = s \bP^x( \tau_{B_1} > s, \xi_\lam > s) \\
 &= s \bP^x( \tau_{B_1} > s) \bP^x(\xi_\lam > s)
 \ge s (1-\eps) e^{-\lam s},
\end{align*}
which yields \eqref{e:GDlb}. \qed

\begin{theorem} \label{T:gdCS}
Suppose $\sX$ satisfies VD and $\UHK(\Psi)$. 
Then $\sX$ satisfies $\FK(\Psi)$ and $\CSA(\Psi)$.
\end{theorem}

\proof
The proof that $\UHK(\Psi)$ plus VD implies $\FK(\Psi)$ is as in
Section 5.5 of \cite{GH}, where the case $\Psi(r)=r^\beta$ is given.

To prove $\CSA(\Psi)$ we will show that $\sX$ satisfies the hypotheses
of Lemma \ref{L:cswk}.
So let $B'=B(x_0, R)$ and $B=B(x_0, R+r)$, and $U=B-B'$,
and let $D_i$, $h$ be as in Lemma \ref{L:hprop}. Set
\begin{align}
 g(x) &= \frac{c_1 h(x)}{\Psi(r)},  \\
 \vp(x) &= 
\begin{cases}
 1 \wedge g(x) & \text{ if } x \in B(x_0, R + r/2)^c, \\
 1  & \text{ if } x \in B(x_0, R + r/2). \\
\end{cases}
\end{align}

Then by Lemma \ref{L:hprop} $\vp =0 $ on $B^c$, and $\vp = 1$ on $B'$,
so it remains to verify the inequality \eqref{e:csa-C}. 

Let $f \in \sF$. Since $g$ is zero outside $U$ we have
\begin{align} \nn 
\int_U f^2 d\Gam(\vp, \vp) &\le \int_U f^2 d\Gam(g, g) 
= \int_\sX f^2 d\Gam(g, g)  \\
\label{e:int1}
&= \int_\sX d\Gam(f^2 g,g) - 2 \int_\sX f g d\Gam(f,g). 
\end{align}
Now writing $\sE_\lam(u,v) = \sE(u,v)+\lam \langle u,v \rangle$,
\begin{align}
\nn
 \int_\sX d\Gam(f^2 g,g)= \sE(f^2 g,g) 
 &\le \sE_\lam(f^2 g,g) \\
\nn
&= c_1 \Psi(r)^{-1} \sE_\lam(f^2 g,G^{D_0}_\lam 1_{D_1}) \\
\label{e:cr1}
&=  c_1 \Psi(r)^{-1} \langle f^2 g,1_{D_1} \rangle  
\le c_1 \Psi(r)^{-1}  \int_U f^2 g dm .
\end{align} 
Here we used \cite[Theorem 4.4.1]{FOT} 
and the fact that $f^2 g \in \sF_{D_0}$ to obtain the third line.
By \eqref{e:cs-lam}, 
\be \label{e:cr2}
 \big| 2 \int_\sX f g d\Gam(f,g) \big| 
 = \half \int_\sX f^2 d\Gam(g,g) + 2 \int_\sX g^2 d\Gam(f,f) .
\ee 
Combining \eqref{e:cr1} and \eqref{e:cr2}, and using the fact that
$g \le c_1$, we obtain
\begin{align*}
 \int_U f^2 \, d\Gamma(g,g) 
&\le 4 \int_U g^2 d\Gamma(f,f) + 2 c_1 \Psi(r)^{-1}  \int_U  g f^2 dm \\
&\le 4c_1^2 \int_U d\Gamma(f,f) + 2 c_1^2 \Psi(r)^{-1}  \int_U f^2 dm.
\end{align*}
Thus the hypotheses of Lemma \ref{L:cswk} hold, and so $\CSA(\Psi)$ 
holds. \qed

\begin{remark} 
{\rm 
While the proof above is based on the argument in 
Section 3 of \cite{BBK}, it is much simpler, since we do
not need to consider the integral over arbitrary balls.
Further, 
the condition $\CS(\Psi)$ requires H\"older continuity of the
cutoff function, and this was proved by  using a parabolic Harnack
inequality, which is equivalent to the full (upper and lower)
heat kernel bounds $\HK(\Psi)$.
It seems unlikely that the conditions VD and $\UHK(\Psi)$ are
sufficient to ensure the existence of a H\"older continuous cutoff
function. 
}
\end{remark}

\ms We conclude
this section by giving a sketch of the proof that $\CSA(\Psi)$ follows from
the condition $\CS(\Psi)$ introduced in \cite{BB3, BBK}.

\begin{lemma} \label{L:cs2csa}
Let $\sX$ satisfy $\VD$. 
Suppose that for every $x \in \sX$ and $r>0$ there
exists a cutoff function $\vp$ for $B(x,r) \subset B(x,2r)$
such that if $f: B=B(x,2r) \to \bR$ then, writing $V=B(x,2r)-B(x,r)$,
\be \label{e:cs-suff}
 \int_V f^2 \, d\Gamma(\vp,\vp) 
\le c_1 \Big( \int_V d\Gamma(f,f) 
+\Psi(r)^{-1} \int_V f^2 dm  \Big).
\ee 
Then $\CSA(\Psi)$ holds.
In particular $\CS(\Psi)$ implies $\CSA(\Psi)$.
\end{lemma}

\proof Let $x_0\in \sX$, $R,r>0$, and 
$B'=B(x_0, R)$ and $B=B(x_0, R+r)$, and $U=B-B'$. 
In view of Lemma \ref{L:cswk} it is enough to prove that there exists 
$c_2< \infty$ such that for $f: U \to \bR$,
\be \label{e:etp5}
 \int_U f^2 \, d\Gamma(\vp,\vp) 
\le c_2 \int_U  d\Gamma(f,f) + c_2 \Psi(r)^{-1} \int_U  f^2 dm . 
\ee

Set $r_0 = r/3$, and let $B(z_i, r_0)$ be a covering of $B'$ by
balls such that $B'_i=B(z_i, r_0/2)$ are disjoint and each $z_i \in B'$. 
Then VD implies there exists $M$ such that any ball $B(y, r_0/100) \subset B$ 
intersects at most $M$ of the balls $B_i$. Let $\vp_i$ be a cutoff function 
for $B'_i \subset B_i = B(z_i, 2 r_0)$ satisfying \eqref{e:cs-suff}. Then
$$ \int_{B_i} f^2 d\Gam(\vp_i,\vp_i) 
\le c_1 \Big( \int_{B_i} d\Gam(f,f)+  \Psi(r_0/3)^{-1} \int_{B_i} f^2 dm \Big). $$
Now set $\vp(x) = \max_i \vp_i(x)$. Then $\vp$ is clearly 1 on $B'$ and 
zero outside $B$.

If $B''=B(y,r_0/100)$ and $B_i$, $i=1, \dots m$ are the balls which intersect
$B''$, then
$$ d\Gam(\vp,\vp) \le \sum_{j=1}^m d\Gam(\vp_j, \vp_j). $$
Thus
\begin{align*}
 \int_{B} f^2 d\Gam(\vp,\vp)  &\le \sum_i \int_{B_i} f^2 d\Gam(\vp_i,\vp_i)  \\
&\le \sum_i  c_1 \Big( \int_{B_i} d\Gam(f,f)+ \Psi(r_0/3)^{-1} \int_{B_i} f^2 dm \Big) \\
&\le c_1 M \Big( \int_{B} d\Gam(f,f)+  \Psi(r_0/3)^{-1} \int_{B} f^2 dm \Big),
\end{align*}
proving \eqref{e:etp5}.\qed

\section{Stochastic Completeness} \label{s:scom}

\noindent{\em Proof of Theorem \ref{T:ScomI}}.
Following Davies \cite[Theorem 7]{D} let  $f\ge 0$ be a function with 
compact support and let $u_t = P_t f$.
We remark that to prove stochastic completeness, by standard 
density arguments it is sufficient to prove that 
\begin{align} \label{e:sc-davies}
 \int_\sX f \, dm \leq \int_\sX u_t \, dm  \q \hbox{ for some $t>0$.} 
\end{align}
Indeed, note that since $P_t$ is  self-adjoint in $L^2(\sX,m)$, this implies
$\langle 1-P_t 1, f \rangle \leq 0$ and therefore $P_t1=1$ m-a.e.

Let $(a_n)$ be an increasing sequence with $a_0=1$, and define
$\vp$, $b_n$, $b^*$ and $C_0$ as in \eqref{e:vpdef0}--\eqref{e:C0def}.
We assume that $(a_n)$ is chosen so that $b^*=1$.
Let $t \in (0,1)$. Then
$$ \langle f, \vp_n \rangle - \langle u_t, \vp_n \rangle 
   = -\int_0^t  \frac{d}{ds} \langle u_s,\vp_n \rangle ds $$
and
$$  - \frac{d}{ds} \langle u_s,\vp_n \rangle = \sE(u_s,\vp_n) 
= \int_\sX d\Gamma(u_s, \vp_n). $$
So, by Cauchy-Schwarz and Proposition \ref{p:CSDG}, and recalling that $t<1$,
\begin{align*} 
 \langle f, \vp_n \rangle - \langle u_t,\vp_n \rangle 
 &= \int_0^t \int_\sX \vp \cdot \vp^{-1}  d\Gamma(u_s, \vp_n) \, ds \\
&\le \Big( \int_0^t \int_\sX \vp^2 \,  d\Gamma(u_s,u_s) \, ds\Big)^{1/2}
    \Big( \int_0^t \int_\sX \vp^{-2} \,  d\Gamma(\vp_n,\vp_n) \, ds \Big)^{1/2} \\
&\le \sqrt{2} \, \| f \vp \|_2 e^{2C_0 t} 
 \big( \sup_{U_n} \vp^{-1}\big) \Big ( \int_\sX d\Gamma(\vp_n,\vp_n) \Big)^{1/2}.
\end{align*}
On $U_n$ we have $a_n \le \vp \le a_{n+1}$, so
$\sup_{U_n} \vp^{-1} \le a_n^{-1}$.
Using $\CSD(D_n, D_{n+1},\th_n)$ with $f=1$,
$$ \int_\sX  d\Gamma(\vp_n,\vp_n) 
= \int_{U_n}   d\Gamma(\vp_n,\vp_n) 
\le \th_n m(U_n).  $$
So,
\be   \label{e:sc-e1}
  \langle f, \vp_n \rangle  - \langle u_t,\vp_n \rangle 
\le \sqrt{2} \, \| f \vp\|_2  
 \exp \Big( 2C_0 t+ \half \log ( \th_n m(U_n)) - \log(a_n) \Big). 
\ee
If there exists a subsequence $(n_k)$ such that
\be \label{e:cssub}
\lim_{k \to \infty} (\langle f, \vp_{n_k} \rangle 
   - \langle u_t,\vp_{n_k} \rangle) \le 0, 
\ee 
then, since 
$$ \int_\sX u_t \, dm= \lim_{k} \int_{\sX} u_t \vp_{n_k} \, dm, $$
we obtain \eqref{e:sc-davies} and so deduce stochastic completeness. 

\noindent (a) If $\th_n \le c_1$ 
we choose $a_n=2^n$, so that $b_n=b^*=1$ and $C_0 =c_1<\infty$.
Then \eqref{e:asybI} implies that the right side of \eqref{e:sc-e1}
converges to 0.

\noindent (b) (Recall in this case that $\th_n = c_0^2 n^2$.) 
Let $\al>0$, and consider sequences
$(a_k)$ such that $C_0 =C_0((a_k))=\al^2$. 
We wish $a_k$ to be as large as possible given these constraints, and
so choose $b_k = 1\wedge (\al/\th_k^{1/2})$. 
Now fix $n \gg 1$, let $m = \lam \log n$ where $\lam>0$, and let 
$\al= c_0 m$. We have
$$ a_n = \prod_{j=1}^n \big(1+ b_j ) \ge 
\prod_{j=m}^n \big( 1+ \frac{m}{j} \big). $$

So since $\log(1+x) \ge \half x$ for $x \in (0,1)$, for $n$ large enough
\begin{align*}
 \log a_n \ge \half m \sum_{j=m}^n j^{-1} 
 \ge \half \lam (\log n) ( \log n - \log \log n -1)
  \ge \fract13 \lam (\log n)^2. 
\end{align*}
Writing $E(n,\lam)$ for the term in the exponential in
\eqref{e:sc-e1}, if $\log m(U_n) \le 2 b (\log n)^2$ then 
\begin{align*}
 E(n,\lam) &=  2 c_0^2 m^2 t + \half \log( \th_n m(U_n)) - \log a_n \\
  &\le  \half \log( \th_n m(U_n)) - (\log n)^2 ( \lam \fract13 - 2c_0^2 \lam^2 t)\\
 &\le  {\log c_0 n} - (\log n)^2 \Big( \fract{\lam}{3} - b -  2c_0^2 \lam^2 t \Big).   
\end{align*}
Choosing $\lam =9b$ and $t$ small enough so that $2c_0^2 \lam^2 t  \le b$, it
follows that 
$$ E(n,\lam) \le -b( \log n)^2 +   \log c_0^2 n , $$
and \eqref{e:cssub} holds. \qed

\begin{remark} \label{R:scgam}
{\rm 
We have just considered the cases $\th_n \le c_1$ and $\th_n = c_1n$, as
for our applications these are of most interest. By arguments similar to the above
it is straightforward to show that if $\th_n = c_0^2 n^{2\gam}$ with $0<\gam <1$,
then stochastic completeness holds provided
\be
 \log m(U_n) \le c(\gamma) n^{1 \wedge (2-2\gam) }. 
\ee

}\end{remark}

\sms 
We now give some examples of the use of the criterion in Theorem \ref{T:ScomI},
and begin by showing that we can recover the result of Davies \cite{D}.

\begin{example}{\rm
Let $\sX$ be a manifold containing a point $0$, and such that there exists
$b>0$ such that 
\be \label{e:sc-dav}
m(B(0,r)) \le e^{b r^2}. 
\ee

Let $(r_n)$ be increasing with $\lim r_n = \infty$.
Set $D_n=B(0,r_n)$ and let $U_n = D_{n+1}-D_n$. Let $\vp_n$ be `linear'
on $U_n$, so that
$$ \vp_n (x) = 1 \wedge \Big( \frac{r_{n+1} - d(0,x) }{r_{n+1}-r_n } \vee 0\Big),
\hbox { and } ||\grad \vp_n ||_\infty = \frac{1}{r_{n+1}-r_n}. $$
Letting $\th_n = (r_{n+1}-r_n)^{-2}$, clearly we have
\be \label{e:csc}
 \int_{U_n} f^2 \, d\Gamma(\vp_n,\vp_n) 
\le  \th_n \int_{U_n} f^2 dm, 
\ee 
and so $\CSD(D_n, D_{n+1}, \th_n)$ holds.
Let $r_n = \log n$, so that $\th_n \sim n^2$.
Then $m(U_n) \le m(D_{n+1}) \le \exp( b (\log(n+1))^2)$, 
so \eqref{e:asyaI} holds and $\sX$ is stochastically complete.
} \end{example}

\begin{remark}
{\rm Improving the condition $\log V(0,r) \le br^2$ to
 $\log V(0,r) \le r^2\log r$ allowed by Theorem \ref{T:vgc} seems to require
more delicate techniques.}
\end{remark}

\section{The pre-Sierpinski carpet}
\label{s:bmsc}

In this section we will give an example of an MMD space which is geodesically
incomplete but stochastically complete. The example is based on the
`pre-Sierpinski carpet' -- see \cite{O1}.

\begin{figure}
\includegraphics{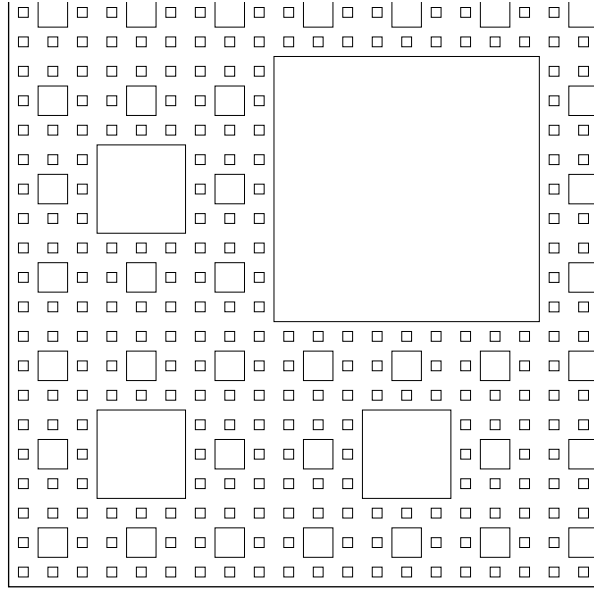} 
\caption{The pre Sierpinski carpet} 
\end{figure}

The standard Sierpinski carpet in $d$ dimensions (with $d \ge 2$) can be 
constructed by an analogue of the construction of the Cantor set.
Starting with $F_0=[0,1]^d$, divide $F_0$ into $3^d$ subcubes each
of side $3^{-1}$, and remove the middle cube; call this set $F_1$.
Repeating this construction, we obtain a decreasing sequence of compact sets
$F_n$; the Sierpinski carpet is defined as 
$$ F = \cap_{n=0}^\infty F_n. $$
Let $M_d = 3^d-1$; then $F$ has Hausdorff dimension
$$ d_f = \frac{ \log M_d}{\log 3}. $$
Note that $F_n$ is a union of $M_d^n$ cubes each of side $3^{-n}$.
Let
\be
   \wt F_n = 3^n F_n = \{ 3^n x: x \in F_n \}, \q 
\wt F = \cup_{n=0}^\infty \wt F_n. 
\ee
The set $\wt F$ is the {\em pre-Sierpinski carpet}, and is a countable
union of copies of the unit cube $[0,1]^d$. The interior of $\wt F$ is a 
standard open domain in $\bR^d$, with a Lipschitz boundary.
We write $\sX= \wt F$, and will take $d \ge 3$.
Let $\mu$ be Lebesgue measure restricted to $\sX$.

We summarise some properties of $\sX$. Let $d(x,y)$ denote the 
shortest path distance in $\sX$. Then (see \cite[Lemma 7.3]{BB1}
for the case $d=2$) we have
\be \label{e:dmetric}
|x-y|\leq d(x,y)\leq c |x-y|, \q x,y \in \sX. 
\ee
We write $B_d(x,r)$ for balls in the metric $d$.
Then (see \cite[Lemma 2.3(e)]{BB2}) we have
\be \label{e:volg}
V_d(x,r)  = \mu(B_d(x,r)) \asymp 
\begin{cases}
  r^d, & 0\le r \le 1, \\
  r^{d_f}, & r>1.
\end{cases}
\ee
In particular $\sX$ satisfies VD.

Now set 
\begin{align*}
\sE(f,f)=\int_\sX |\nabla f|^2 \, d\mu, \qquad f\in H^1(\sX),
\end{align*}
where $H^1=H^1(\sX)$ denotes the set of functions $f$ for which 
$\int_\sX f^2 \, dx+\sE(f,f)<\infty$. 
Then $(\sE, H^1)$ is a regular local Dirichlet form on $L^2(\sX,\mu)$.
The associated
Hunt process $W=(W_t, t\ge 0, \bP^x, x \in \sX)$ is Brownian motion
in $\sX$ with normal reflection on the boundary $\pd \sX$.
For the existence and uniqueness
in law of this process we refer to \cite{BH}. The process is
reversible with respect to $\mu$, and its 
generator is given by the Neumann Laplacian $\Delta$ on $\sX$.

Let $p_t(x,y)$ denote the heat kernel associated with $W$.
Many of the properties of $W$ and $p_t$ can be summarised by two indices. 
The first is $d_f$, the Hausdorff dimension of the space $F$. The second,
denoted $d_w$, and called the {\em walk dimension}, gives the long range
space-time scaling on $\sX$. For Sierpinski carpets in $d\ge 3$  this satisfies 
 $2< d_w < d_f$ -- see \cite[Section 5]{BB2}.
Let $\Psi = \Psi_{2, d_w}$ be as defined in \eqref{e:Psidef}.

\begin{theorem} \label{T:Wprop}
(a) $(\sX, \sE)$ satisfies $\HK(\Psi)$. \\
(b) $W$ has a Greens function $g(x,y)$ such that there
exist positive constants $c_1$-$c_4$ such that
\begin{align*}
c_1 |x-y|^{2-d} & \leq g(x,y)  \leq c_2  |x-y|^{2-d}\quad \text{if $|x-y|\leq 1$}, \\
c_3 |x-y|^{d_w-d_f} & \leq g(x,y)  \leq c_4  |x-y|^{d_w-d_f}\quad \text{if $|x-y|> 1$} .
\end{align*}
(c) The conditions $\CS(\Psi)$ and $\CSA(\Psi)$ hold for $(\sX,\sE)$.
\end{theorem}

\proof (a) is proved in \cite[Theorem 6.9]{BB2}, and (b) in
\cite[Corollary 6.10]{BB2}. That $\CS(\Psi)$ holds follows from
\cite{BBK}. $\CSA(\Psi)$ then follows by Lemma \ref{L:cs2csa}, or
alternatively by Theorem \ref{T:main}. \qed

\sms
Let $a(x)>0$, $x\in \sX$ be a real-valued function on $\sX$. Then, 
we define the additive functional 
\[
A_t=\int_0^t \frac 1 {a(W_s)} \, ds
\]
and the time-changed process $Y=Y^{(a)}$ by
\[
Y_t=W_{\tau_t}, \qquad t\geq 0,
\]
where $(\tau_t)$ denotes the inverse of $(A_t)$. The process $Y$ is symmetric 
with reversible measure $m(dx)=m_a(dx)=a^{-1}(x)\, \mu(dx)$ and its generator 
$\sL_a$ satisfies
$$ \int g \sL_a f  a^{-1} \, d\mu 
=  \langle \sL_a f,g\rangle_{L^2(m_a)} = -\sE(f,g) = -\int \grad f \cdot \grad g  d\mu
 =  \int (\Delta f ) g d\mu , $$
so that
\be 
 \sL_a f = a \Delta f.
\ee
The Dirichlet form associated with $Y$ is the form $(\sE, \sD_a)$ 
on the base space $L^2(\sX,m_a)$. 
Here $\sD_a$ is the closure of $C^1_0(\sX)$ with respect to 
$N_a(f) = \sE(f,f) + ||f||^2_{L^2(m_a)}$.
We refer to this form as $\sE_a$ for short. 
Recall from \eqref{e:intrin} the definition of the 
intrinsic metric $\varrho_a$ associated with $\sE_a$; we have
\be
 \rho_a(x,y)=\sup \{ u(x)-u(y):\, u\in \sM_a, \},
\ee
where
\be
\sM_a = \{ u \in \sD_a \cap C(\sX):
\, |\nabla u|^2 \leq a^{-1} \}.
\ee

\sm Let $p>0$. We now just consider the case 
\be
  a(x) = 1 \vee d(0,x)^p.
\ee

The main result of this section is the following. Recall that
we have $d\ge 3$, and that $2 < d_w < d_f$.

\begin{theorem} \label{thm_sc}
(i) The process $Y=Y^{(p)}$ is stochastically complete if and only 
if $p\leq d_w$. \\ 
\noindent (ii) On the other hand, (VGC) holds if and only if 
$p\leq 2$ or $p>d_f$. 
In particular, for $p\in (2,d_w)$ the process $Y$ is stochastically 
complete but (VGC) fails.
\end{theorem}

We begin by relating the metrics $\rho_a$ and $d$ on $\sX$.

\begin{lemma} \label{L:rhoest}
Let $R>0$, and $x \in \pd B_d(0,R)$, $y \in \pd B_d(0, 2R)$.
Then
\be
   \rho_a(x,y) \asymp  R^{1- p/2} , \hbox{ if } R\ge 1,
\ee
while $\rho_a(x,y) \asymp R$ if $R\in [0,1]$.
\end{lemma}

\proof By Theorem 4.1 in Chapter 5 of \cite{St2} we have
\be
\varrho_a(0,x) = \inf_\gamma \int_0^1 | \overset{.}{\gamma}|  \, 
\frac 1 {\sqrt{a(\gamma(s))}} \, ds.
\ee
If $\gam$ is any path in $B_d(0,2R)-B_d(0,R)$ then
$$ \int_0^1 | \overset{.}{\gamma}|  \, 
\frac 1 {\sqrt{a(\gamma(s))}} \, ds \asymp R^{-p/2} \int \dot \gam
= R^{-p/2} |\gam|, $$
where $|\gam|$ denotes the length of $\gam$.
It follows that $\rho_a(x,y) \ge c R^{1-p/2}$.

For the upper bound, the geometry of the pre-carpet implies that
if $C>1$ is large enough then
we can find a path $\gam_1$ between $x$ and $y$ which lies inside
$B_d(0,CR)-B_d(0,R)$ and has length less than $c_1 R$. Therefore
$\rho_a(x,y) \le c_2 R^{1-p/2}$.\qed

\begin{proposition} \label{p:metr-meas} 
The metric $\rho_a$ and measure $m_a$ satisfy the following. \\
(i) $\varrho_a(0,\infty)=\infty$ if and only if $p\leq 2$. 
In particular, $(\sX, \rho_a)$ is not geodesically complete when $p>2$. \\
(ii) $m_a(\sX)=\infty$  if and only if $p\leq d_f$.
\end{proposition}

\proof
(i) Let $x_k=(2^k, 0, \dots, 0)$ be the points on $\pd B_d(0, 2^k)$. 
Then by Lemma \ref{L:rhoest} we have $\rho_a(x_k, x_{k+1}) \asymp 2^{k(1-p/2)}$.
If $p \le 2$ the sum $\sum_k \rho_a(x_k, x_{k+1})$  diverges, and hence $\rho_a(0, x_k ) \to \infty$,
while if $p>2$ then $\lim_k \rho_a(0, x_k ) \le C_1 < \infty$; (i) then
follows.

\noindent (ii) By \eqref{e:volg}
\begin{align*}
m_a(\sX) &= \sum_k m_a(B_d(0, 2^{k+1})- B_d(0, 2^k))  \\
 &\asymp \sum_k 2^{-kp} \mu(B_d(0, 2^{k+1})- B_d(0, 2^k))
  \asymp \sum_k 2^{-kp} 2^{kd_f},
\end{align*}
which is infinite if and only if $p\leq d_f$. \qed

\sms
We now look at (VGC) for the metric measure space $(\sX, \rho_a, m_a)$.
We set $h^{-1}(r):=\int_0^r (1\vee s)^{-p/2}  \, ds$, $r\geq 0$, so 
that $\varrho_a(0,x)\asymp h^{-1}(|x|)$ if $|x|\geq 1$. Further, let 
$h(r):=\inf\{ s: \, h^{-1}(s)> r\}$ be the right-continuous inverse. In particular,
\[
B_d(0,h(c_2r))\subseteq B_{\rho_a}(0,r) \subseteq B_d(0,h(c_1r)), \qquad  r>1. 
\]
Moreover, if $p\leq 2$ we have $\lim_{r\to \infty}h^{-1}(r)=\infty$, 
so $h(r)<\infty$ for all $r$. 
On the other hand,
if $p>2$ we have that $R_0:=\lim_{r\to \infty} h^{-1}(r)<\infty$, thus $h(r)=\infty$ 
and $B_{\rho_a}(0,r)=F$ for $r\geq R_0$.

\begin{lemma} \label{l:est_maballs}
Let $p\leq 2<d_f$.There exist positive constants $c_1$--$c_5$ and  
$r_0$ such that for all $r>r_0$ 
\[
c_1 h(c_2 r)^{d_f-p} \leq m_a(B_{\rho_a}(0,r))\leq  c_3 h(c_4 r)^{c_5(d_f-p)}.
\]
\end{lemma} 

\proof
The lower bound is immediate from \eqref{e:volg} as
\begin{align*}
m_a(B_{\rho_a}(0,r))
 &\geq \int_{B_d(0,h(c r))} \frac 1 {a(x)} \mu(dx) \\
 &\geq \mu(B_d(0,1))+ h(r)^{-p} 
     \left(\mu(B_d(0,h(c r)))-\mu(B_d(0,1))\right)\\ 
&\geq c+c h(cr)^{d_f-p}
\end{align*}
for $r$ sufficiently large, where we used the fact that $h$ is increasing.

To prove the upper bound note that for  $k\ge 0$ we have 
$a(x) \ge 2^{kp}$ on the set 
$U_k=B_d(0,2^{k+1}) - B_d(0,2^k)$. 
Let $k_0(r) =\min\{k: 2^k \ge h(cr) \}$. 
Then, for all $r$ large enough we have again by \eqref{e:volg}
\begin{align*}
m_a(B_{\rho_a}(0,r))&\leq\int_{B_d(0,h(c r))} \frac 1 {a(x)} \mu(dx) \leq  \mu(B_d(0,1))
 +  c \sum_{k=0}^{k_0(r)} m_a(U_k) \\
&\le c+ c \sum_{k=0}^{k_0(r)} 2^{k(d_f-p)} 
\le c h(c r)^{c(d_f-p)}.
\end{align*}
\qed

\sm {\em Proof of Theorem \ref{thm_sc}}\\
(i) First let $p\le d_w$.
Let $R>1$, $R_n=R^n$, $D_n := B_d(0, R_n)$ and $U_n=D_{n+1}-D_n$.  
Thus 
$$ m_a(U_n) \asymp R_n^{-p} \mu(U_n) \asymp R_n^{ d_f -p}. $$
By Theorem \ref{T:Wprop} we have $\CSA(\Psi)$ for the space $(\sX,d,\sE,\mu)$.
So there exists a cutoff function $\vp_n$ for $D_n \subset D_{n+1}$
such that if $f: U_n \to \bR$ then
\begin{align*}
 \int_{U_n} f^2 d\Gam(\vp_n,\vp_n) 
&\le \frac18 \int_{U_n} \vp_n^2 d\Gam(f,f) + c_0 R_n^{-d_w} \int_{U_n} f^2 d\mu \\
&\le \frac18 \int_{U_n} \vp_n^2 d\Gam(f,f) + c_1 R_n^{p-d_w} \int_{U_n} f^2 dm_a. 
\end{align*}
Thus in the space $(\sX, \rho_a, m_a)$, $\CSD( D_n, D_{n+1}, \th_n)$ holds with
$\th_n = c_1 R_n^{p-d_w}$.
As $p\le d_w$ we have $\th_n \le c_1$ and hence by Theorem \ref{T:ScomI}(a)
stochastic completeness holds provided \eqref{e:asybI} holds. However,
$$ \frac{\th_n m_a(U_n)}{4^n } \le c_1 4^{-n} R_n^{-d_w +p} R_n^{d_f -p}
 = c_1 (R^{d_f-d_w}/4)^n, $$
and taking $R$ small enough so that $R^{d_f-d_w}<4$ it follows that
stochastic completeness holds. 

\sms
Now we consider the case $p>d_w$. Since the process $W$ is 
stochastically complete, from the definition 
of stochastic completeness it is immediate that $Y$ is stochastically complete 
if and only if  $A_\infty=\infty$ $\bP^x$-a.s. for any $x$.  
Note that $A_\infty=\infty$ is a tail event, i.e.\ it 
is in $\sigma(W_s, s\geq t)$ for all $t$, so 
$\bP^x[A_\infty=\infty]$ is either 0 or 1 for 
all $x$ (cf.\ Theorem 8.7 in \cite{BB2}).

Let $D_0 =  B_d(0,1)$ and for $n \ge 1$ set
$D_n = B_d(0, 2^n)-B_d(0, 2^{n-1})$. Then using
the bounds for the Green kernel $g(x,y)$ of $W$ in 
Theorem \ref{T:Wprop}(b), 
\begin{align*}
\bE^0 A_\infty &= \int_\sX a(x)^{-1} g(0,x) \mu(dx)  \\
&\le c \int_{B_d(0,1)} |x|^{2-d} dx 
 + c \sum_{n=1}^\infty \int_{D_n} d(0,x)^{-p+d_w-d_f} dx  \\
& \leq  c + c \sum_{n=1}^\infty 2^{ n(-p+d_w-d_f)} 2^{n d_f} 
\leq c+ c \sum_{n=1}^\infty  2^{ n(d_w-p)} < \infty, 
\end{align*}
Hence, $A_\infty<\infty$ $\bP^0$-a.s., 
and so $Y$ is stochastically incomplete.

\sm (ii) 
Let us first consider the case $p<2$. 
Then, $d_f-p>0$ and $h(r)\asymp r^{\gam}$ for large $r$ 
with $\gam:=(1-p/2)^{-1}$. Hence, we use \eqref{e:volg} to obtain
\begin{align*}
\int_1^\infty \frac r {\log m_a(B_{\rho_a}(0,r))} dr 
\geq \int_1^\infty \frac r {c_1+ c_2 \log r} \, dr=\infty.
\end{align*}
If $p=2$, $h(r)=e^{r-1}$, so
\begin{align*}
\int_1^\infty \frac r {\log m_a(B_{\rho_a}(0,r))} dr 
\geq \int_1^\infty \frac r {c_1+ c_2 (r-1)} \, dr=\infty.
\end{align*}
Finally, in the case $p>2$ we have that $R_0:=\lim_{r\to \infty} h^{-1}(r)<\infty$, 
thus $h(r)=\infty$ and $B_{\rho_a}(0,r)=\sX$ for $r\ge R_0$. In particular, by 
Proposition \ref{p:metr-meas} ii) we get for such $r$ that  
$m_a(B_{\rho_a}(0,r))=m_a(\sX)=\infty$ if and only $p\leq  d_f$. Hence,
\begin{align*}
\int_1^\infty \frac r {\log m_a(B_{\rho_a}(0,r))} dr 
 \begin{cases}
 < \infty & \text{if $p\leq  d_f$,} \\
 =\infty & \text{if $p> d_f$.}
 \end{cases}
\end{align*}
\qed

\sm SA: Institut f\"ur Angewandte Mathematik \\
Rheinische Friedrich-Wilhelms-Universit\"at Bonn \\
Endenicher Allee 60, 53115 Bonn, Germany. \\
andres@iam.uni-bonn.de

\sm MB: Department of Mathematics,
University of British Columbia, \\
Vancouver, B.C. V6T 1Z2, Canada. \\
barlow@math.ubc.ca

\end{document}